\documentclass[onefignum,onetabnum]{siamart250211}

\usepackage{lipsum}
\usepackage{amsfonts}
\usepackage{graphicx}
\usepackage{epstopdf}
\usepackage{algorithmic}
\usepackage{comment}
\usepackage{amssymb}
\usepackage{subcaption}

\ifpdf
\DeclareGraphicsExtensions{.eps,.pdf,.png,.jpg}
\else
\DeclareGraphicsExtensions{.eps}
\fi

\newsiamremark{remark}{Remark}
\newsiamremark{hypothesis}{Hypothesis}
\crefname{hypothesis}{Hypothesis}{Hypotheses}
\newsiamthm{claim}{Claim}

\DeclareSymbolFont{lettersA}{U}{txmia}{m}{it}
\newcommand{\R}{\mathbb{R}}

\renewcommand{\(}{\left(}
\renewcommand{\)}{\right)}
\newcommand{\yvec}{\mathbf{y}}
\newcommand{\fvec}{\mathbf{f}}
\newcommand{\zvec}{\mathbf{z}}
\newcommand{\ones}{\mathbf{1}}

\definecolor{codegreen}{rgb}{0,0.6,0}
\definecolor{codegray}{rgb}{0.5,0.5,0.5}
\definecolor{codepurple}{rgb}{0.58,0,0.82}
\definecolor{backcolour}{rgb}{0.95,0.95,0.92}

\headers{Extended STS Methods}{D. R. Reynolds, S. Amihere, and M. Aggul}

\title{Implicit-Explicit and Split-Explicit Super-Time-Stepping Methods\thanks{Submitted to the editors \date.
\funding{This work was supported by the U.S. Department of Energy, Office of Science, Office of Advanced Scientific Computing Research, Scientific Discovery through Advanced Computing (SciDAC) Program through the Computational Evaluation and Design of Actuators for Core-Edge Integration (CEDA) project, and through the Frameworks, Algorithms and Software Technologies for Mathematics (FASTMath) Institute.}}}

\author{Daniel R. Reynolds\thanks{University of Maryland Baltimore County, Department of Mathematics \& Statistics, Baltimore, MD, USA (\email{dreynolds@umbc.edu}, \email{samihere@umbc.edu})}
  \and Sylvia Amihere\footnotemark[2]
  \and Mustafa Aggul\thanks{Southern Methodist University, Department of Mathematics, Dallas, TX, USA (\email{maggul@smu.edu})}
}

\usepackage{amsopn}

\makeatletter
\newcommand*{\addFileDependency}[1]{
  \typeout{(#1)}
  \@addtofilelist{#1}
  \IfFileExists{#1}{}{\typeout{No file #1.}}
}
\makeatother

\ifpdf
\hypersetup{
  pdftitle={Implicit-Explicit and Split-Explicit Super-Time-Stepping Methods},
  pdfauthor={D. R. Reynolds, S. Amihere, and M. Aggul}
}
\fi

\begin{document}

\maketitle

\begin{abstract}
  Multiphysics initial-value problems couple processes with distinct stability properties, such as advection, diffusion, and stiff local reactions.  Standard implicit-explicit (ImEx) additive Runge--Kutta (ARK) methods can treat these processes accurately, but require globally coupled implicit solves when diffusion is grouped with reaction; operator splitting avoids such solves but typically provides weaker coupling and no inexpensive temporal error estimate; and PIROCK is tied to a specific Runge--Kutta--Chebyshev super-time-stepping (STS) construction.  We introduce \emph{extended super-time-stepping} (ExtSTS) methods, a family of time integration schemes that combine super-time-stepping methods for diffusive terms with explicit, implicit, or ImEx Runge--Kutta treatment of the remaining terms.  The coupling is based on multirate infinitesimal techniques, yielding solve-decoupled methods that retain localized implicit solves, support embedded error estimation for adaptive time stepping, and allow flexible use of modern STS methods.  We present the ExtSTS method family, provide a robust technique for ExtSTS method creation, formulate the corresponding linear stability theory, and construct embedded ImEx, explicit, and implicit ExtSTS methods.  Numerical experiments on one- and two-dimensional advection-diffusion-reaction problems show that ExtSTS methods are robust across parameter regimes and boundary conditions, and are often more efficient than ARK, Strang splitting, and PIROCK methods, especially when strong coupling between operators is important.
\end{abstract}



\section{Introduction}
Recent years have seen the development of numerous methods for systems of ordinary differential equation (ODE) initial-value problems (IVP) that are designed to exploit underlying structure in the dynamical system.  Beginning with an IVP of the form
\begin{equation}
  \label{eq:IVP_unsplit}
  \yvec'(t) = \fvec(t,\yvec), \quad \yvec(t_0) = \yvec_0,
\end{equation}
these methods typically assume that the right-hand side function $\fvec$ may be additively split into a set of functions, $\fvec(t,\yvec) = \sum_{i} \fvec^{i}(t,\yvec)$, such some terms $\fvec^i$ may be stiff with others nonstiff, or some may incorporate dynamics that naturally evolve on a time scale that is faster or slower than others.  Such structure-exploiting methods include \emph{additive Runge--Kutta} (ARK) methods \cite{Ascher1997} and exponential methods \cite{hochbruck_exponential_2010} that evolve all $\fvec^i$ using a shared step size, but where stiff and nonstiff components may be treated separately.  \emph{Multirate} methods \cite{gear_multirate_1984,knoth_implicit-explicit_1998}, on the other hand, evolve $\fvec^i$ using different time steps, thus allowing more rapidly-evolving processes to be evolved using smaller time steps than slower-evolving dynamics.

In this work, we consider IVPs with up to three partitions,
\begin{equation}
  \label{eq:IVP_3split}
  \yvec'(t) = \fvec^{A}(t,\yvec)+\fvec^{D}(t,\yvec)+\fvec^{R}(t,\yvec), \quad \yvec(t_0) = \yvec_0.
\end{equation}
Here, we focus on algorithms that evolve all partitions using a shared step size, $h$, but that treat each $\fvec^i$ using a different type of time integration method. Specifically, we assume that $\fvec^D$ encapsulates diffusive processes such that the spectrum of its Jacobian is located on (or near) the negative real axis, i.e. $\lambda(\partial_{\yvec} \fvec^D) \subset \R^{-}$, $\fvec^A$ corresponds with advective transport processes, and $\fvec^R$ corresponds with stiff reactive processes.  When encountering each of these operators in isolation, it is typical to apply implicit or stabilized explicit \emph{Super Time Stepping} (STS) Runge---Kutta methods \cite{van_der_houwen_internal_1980,verwerRKCTimesteppingAdvection2004}, explicit Runge--Kutta methods \cite{Hairer1993}, and implicit $L$-stable methods \cite{hairer_solving_1996}, respectively.  Our proposed \emph{extended super time stepping} (ExtSTS) family of methods do just that, but in a manner that ensures high order accurate coupling between these processes, as well as a natural temporal error estimator to enable step size adaptivity.

The methods we propose here are not the first to examine problems of the form \eqref{eq:IVP_3split}.  The most direct competitors are the previously-mentioned ARK methods, \emph{operator splitting} methods \cite{Strang1968,Marchuk1968}, and the PIROCK method by Abdulle and Vilmart \cite{abdullePIROCKSwissknifePartitioned2013a}.  Depending on how the operator-splitting method is set up, each of these evolve $\fvec^A$ using an explicit Runge--Kutta method, and evolve $\fvec^R$ using a diagonally implicit Runge--Kutta method; their primary differences lie in how they evolve $\fvec^D$ and how they couple the processes together. In the numerical tests in \cref{sec:results} we compare the performance of our proposed methods against each of these, and provide a detailed explanation of how each method is applied to our benchmark problems.  However, to motivate our goals in proposing ExtSTS methods, we first note flaws in each of these methods that we wish for ExtSTS methods to overcome.

When applying ARK methods to the problem \eqref{eq:IVP_3split}, it is typical to group all stiff operators together to be treated implicitly, $\fvec^R+\fvec^D$, leaving the nonstiff operator $\fvec^A$ to be evolved explicitly.  However, for applications where $\fvec^R$ involves only spatially-local terms, the combination $\fvec^R+\fvec^D$ requires the solution of implicit algebraic systems of equations that are coupled across the entire computational domain.  Modern high-performance computing architectures rely on GPU accelerators that perform excellently on arithmetically-intense local calculations, but that suffer when calculations are frequently synchronized or require global communication.  Thus, we would prefer time integration methods that primarily rely on decoupled nonlinear algebraic solvers and point-to-point communication, in lieu of global implicit solves.

Operator splitting methods, on the other hand, allow perfect flexibility to treat each operator in isolation using whichever algorithm is desired; in our case that would correspond with STS methods for $\fvec^D$ due to their lack of a global implicit solve, explicit methods for $\fvec^A$, and $L$-stable diagonally implicit Runge--Kutta methods for $\fvec^R$.  However, these methods do not admit inexpensive estimates of their temporal error, complicating time step adaptivity.  Furthermore, due to their weak coupling between operators, it is our experience that even higher-order operator splitting methods suffer from large error constants (as will be seen in our numerical results), and thus methods that allow stronger coupling between operators may be desired.

The second-order PIROCK method, on the other hand, addresses both the aforementioned problems of ARK and operator splitting methods, as they effectively apply STS methods to $\fvec^D$, follow that up with an ARK time step for $\fvec^A$ and $\fvec^R$, and include an embedding for estimating the temporal error.  However, these are rooted in the ROCK2 variant \cite{abdulleSecondOrderChebyshev2001} of Runge--Kutta--Chebyshev (RKC) methods, and as such do not allow use of more recent STS algorithms such as Runge--Kutta--Legendre (RKL) \cite{meyerStabilizedRungeKutta2014a} or Runge--Kutta--Gegenbauer (RKG) \cite{osullivanRungeKuttaGegenbauer2019,skarasSupertimesteppingSchemesParabolic2021,tan_explicit_2025}, and their implementation relies on large coefficient lookup tables.  Furthermore, although PIROCK increases the level of coupling over operator splitting methods, be believe that this level of coupling may still be insufficient for applications where the operators $\fvec^A$, $\fvec^D$ and $\fvec^R$ may oppose one another, as will be shown in our results that follow. We conjecture it is for these reasons that even though PIROCK was introduced over a decade ago, it has not seen widespread adoption in the literature.  Thus, we believe that methods which have a similar goal as PIROCK but that allow both a greater flexibility of STS methods, stronger inter-operator coupling, and a more straightforward implementation will allow for more widespread adoption among the computational science community.

Our goal in the present work is to define a family of methods that attain these goals.  We build these by combining STS methods with ARK methods, where instead of coupling them via standard partitioned Runge--Kutta theory, we instead couple them using techniques derived for \emph{multirate infinitesimal} (MRI) methods.  Thus, in the next two subsections \ref{sec:STS-methods} and \ref{sec:MRI-methods}, we highlight the relevant aspects of STS and MRI methods.  In section \ref{sec:PIROCK} we then present the structure of a PIROCK time step for the IVP \eqref{eq:IVP_3split}.  In section \ref{sec:extsts} we propose the new ExtSTS methods, and analyze its accuracy and linear stability in section \ref{sec:analysis}.  We introduce coefficients for candidate implicit-explicit, explicit, and implicit ExtSTS methods in section \ref{sec:methods}.  We then compare the performance of ExtSTS methods against competing approaches in section \ref{sec:results}, and conclude in section \ref{sec:conclusions}.

We note that a class of \emph{multirate Runge--Kutta--Chebyshev} (mRKC) methods were recently introduced in \cite{abdulle_explicit_2022}.  While on the surface these seem related to ExtSTS methods, in that they combine the ideas of multirate and STS methods, they are very different.  mRKC methods assume all problem components are diffusive (i.e., their Jacobians have negative real-valued eigenvalues) and merely combine one RKC method inside another, through applying a single RKC time step (the ``fast'' scale) for the full problem, to construct a homogenized right-hand side evaluation for each stage of a ``slow'' RKC method.  Thus, unlike ExtSTS methods, mRKC methods do not support implicit solvers for any component, and do not support applications with advective terms.

\subsection{STS Methods}
\label{sec:STS-methods}
STS methods were originally invented for diffusion problems of the form \eqref{eq:IVP_unsplit}, where the Jacobian of the right-hand side function satisfies $\lambda(\partial_{\yvec}\fvec)\subset \R^{-}$.  In lieu of applying implicit methods, STS methods add stages to enlarge the stability of explicit methods along $\R^{-}$, thereby preserving the simplicity and locality of explicit time stepping while delivering enhanced stability for diffusion-dominated regimes.  We consider the second-order RKC \cite{sommeijer_rkc_1998} and the RKL \cite{meyerSecondorderAccurateSuper2012,meyerStabilizedRungeKutta2014a} STS methods in this work, both of which have the form
\begin{subequations}
  \label{eq:RKC2}
  \begin{align}
    \label{eq:RKC2a}
    \zvec_0 &= \yvec_n,\\
    \label{eq:RKC2b}
    \zvec_1 &= \zvec_0 + h \tilde{\mu}_1 \fvec(t_n, \zvec_0),\\
    \label{eq:RKC2c}
    \zvec_j &= \mu_j \zvec_{j-1} + \nu_j \zvec_{j-2} + (1-\mu_j-\nu_j)\zvec_0 \\
    \notag
    \label{eq:RKC2d}
    &\quad + h\tilde{\mu}_j \fvec(T_{j-1}, \zvec_{j-1}) + h\tilde{\gamma}_j \fvec(t_{n}, \zvec_0), \quad j=2,\ldots,s,\\
    \yvec_{n+1} &= \zvec_s,
  \end{align}
\end{subequations}
where $T_j = t_n + h c_j$, and the parameters define the specific STS method; analytical expressions can be found in \cite{sommeijer_rkc_1998, meyerSecondorderAccurateSuper2012}.  We note that the more recently introduced second-order RKG method \cite{skarasSupertimesteppingSchemesParabolic2021} has this same structure, and could easily be inserted into our proposed framework.  For both RKC and RKL methods, the extent of the stability region grows proportionally to $s^2$, allowing them to take stable time steps that are substantially larger than those permitted by standard explicit Runge--Kutta schemes in diffusion-dominated problems.  Furthermore, we note that for a given $s$, the method \eqref{eq:RKC2} may be represented as a standard explicit Runge--Kutta method with Butcher tableau $\{A,b,c\} \in \R^{s\times s}\times \R^{s} \times \R^{s}$, and thus their linear stability function for the standard Dahlquist test problem $y'(t) = \lambda y,$ $y_0=1$, with $\lambda \in \R^{-}$ and time step $h>0$ combined as $\eta := h \lambda$, is given by
\begin{equation}
  \label{STS:stability}
  \mathcal{R}_{STS}(\eta) = 1 + \eta b^T(I - \eta A)^{-1}\ones.
\end{equation}

\subsection{The PIROCK Method}
\label{sec:PIROCK}

The method that most closely resembles those proposed here is PIROCK \cite{abdullePIROCKSwissknifePartitioned2013a}, introduced by Abdulle and Vilmart in 2013.  While we do not extend PIROCK methods, we provide their structure here to facilitate comparison later.  For problems of the form \eqref{eq:IVP_3split}, PIROCK starts by taking an initial set of $s-2$ explicit stages of the form \eqref{eq:RKC2a}-\eqref{eq:RKC2c} for only the $\fvec^D$ component.  Following these, PIROCK employs a ``finishing procedure'' to include $\fvec^R$ and $\fvec^A$,
\begin{equation}
  \label{eq:PIROCK}
  \begin{split}
    \zvec_{s-1} &= \zvec_{s-2} + \sigma_{\alpha}h \fvec^D_{s-2},\\
    \zvec_{s} &= \zvec_{s-1} + \sigma_{\alpha}h \fvec^D_{s-1},\\
    \zvec_{s+1} &= \zvec^* + \gamma h \fvec^R_{s+1},\\
    \zvec_{s+2} &= \zvec^* + \beta h \fvec^D_{s+1} + h\fvec^A_{s+1} + (1-2\gamma) h \fvec^R_{s+1} + \gamma h \fvec^R_{s+2},\\
    \zvec_{s+3} &= \zvec^* + (1-2\gamma)h\fvec^A_{s+1} + (1-\gamma) h \fvec^R_{s+1},\\
    \zvec_{s+4} &= \zvec^* + \frac13 h \fvec^A_{s+1},\\
    \zvec_{s+5} &= \zvec^* + \frac{2h\beta}{3} \fvec^D_{s+1} + \frac{2h}{3} \(J^R_s\)^{-1}\fvec^A_{s+4} + \(\frac23 - \gamma\)h\fvec^R_{s+1} + \frac{2h\gamma}{3}\fvec^R_{s+2},\\
    \yvec_{n+1} &= \zvec_s - h\sigma_{\alpha}\(1-\frac{\tau_{\alpha}}{\sigma_{\alpha}^2}\)\(\fvec^D_{s-1}-\fvec^D_{s-2}\) + \frac{h}{4}\fvec^A_{s+1} + \frac{3h}{4}\fvec^A_{s+5} +  \frac{h}{2}\fvec^R_{s+1} + \frac{h}{2}\fvec^R_{s+2} \\
    \notag &\quad + \frac{h}{2-4\gamma}\(J^R_s\)^{-\ell}\(\fvec^D_{s+3}-\fvec^D_{s+1}\),
  \end{split}
\end{equation}
where $J^R_s := I-h\gamma \fvec^R_{\yvec}(\zvec^*)$, and we denote each evaluation of $\fvec^*$ using a subscript corresponding to its $(t,\zvec)$ argument, e.g., $\fvec^A_{j} := \fvec^A(t_n+c_{j}h,\zvec_{j})$.

\subsection{Multirate Infinitesimal Methods}
\label{sec:MRI-methods}

The theoretical foundation of our proposed methods relies on recent work in multirate infinitesimal methods.  The first of this class were \emph{Multirate Infinitesimal Step} (MIS) methods, developed by Wensch, Knoth and Galant, that were constructed as a generalization of split-explicit methods in the context of numerical weather prediction \cite{Wensch2009,Klemp1978,Skamarock1992}. A key contribution from \cite{Wensch2009} was the development of a systematic approach to the order conditions for split-explicit methods based on partitioned Runge--Kutta theory \cite{Wensch2009}, allowing the development of second and third order MIS methods for a variety of applications \cite{Schlegel2012,Schlegel2012a,Schlegel2012b,Knoth2012,Knoth2014,knoth_implicit-explicit_1998,Schlegel2009,Wensch2009}.  The basic approach of these methods is to consider a splitting of \eqref{eq:IVP_unsplit} into slow and fast components, $\mathbf{y}'(t) = \fvec^{F}(t, \mathbf{y}) + \fvec^{S}(t, \mathbf{y})$, where $\fvec^{F}$ represents the fast components and $\fvec^{S}$ represents the slow components.  When evolving a time step $y_n\approx y(t_n) \to y_{n+1} \approx y(t_n + h)$, these methods evaluate $\fvec^S$ infrequently, storing the results to construct forcing terms for a sequence of modified IVPs at the fast time scale,
\begin{subequations}
  \label{eq:mis}
  \begin{align}
    &\notag\text{Let: } z_1 := y_n \\
    &\notag \text{For } i = 2,\ldots,s:\\
    &\label{eq:mis_fast_ode} \quad
    \begin{cases}
      \text{Let: } T_{i} := t_n + c_i h_n, \Delta c_i := c_i - c_{i-1}, \text{ and } v_i(T_{i-1}) := z_{i-1},\\
      \text{Define: } r_i = \frac{1}{\Delta c_i}\sum\limits_{j=1}^{i-1} \(a_{i,j}-a_{i-1,j}\) \fvec^S_j.\\
      \text{Evolve: } v_i'(t) = \fvec^F \left(t,\, v_i\right) + r_i \text{ for } t \in [T_{i-1}, T_i]
    \end{cases}\\
    &\notag \text{Let}: y_{n+1} := z_{s}.
  \end{align}
\end{subequations}
where $(A,b,c) \in \R^{s \times s} \times \R^s \times \R^s$ corresponds to a ``slow'' Butcher tableau that satisfies the conditions:
\begin{align}
  \label{eq:internally_consistent}
  &\sum_{j=1}^{s} a_{i,j} = c_i, \quad i = 1,\ldots,s && \text{[internally consistent]},\\
  \label{eq:explicit}
  &a_{i,j} = 0 \text{ for } i < j && \text{[explicit]},\\
  \label{eq:sorted_abscissae}
  &0 = c_1 < c_2 < \cdots < c_s = 1 && \text{[sorted abscissae]},\\
  \label{eq:stiffly_accurate}
  &a_{s,j} = b_j, \quad j = 1,\ldots,s && \text{[stiffly accurate]}.
\end{align}
If the slow Butcher table includes embedding coefficients $d\in\R^s$, then a stage $s+1$ may be added with the forcing term
\begin{equation}
  \label{eq:MIS_embedding_forcing}
  r_{s+1} = \frac{1}{\Delta c_s}\sum\limits_{j=1}^{s-1} \(d_{j}-a_{s-1,j}\) \fvec^S_j,
\end{equation}
and the corresponding stage provides the embedded solution, $\tilde{y}_{n+1} = z_{s+1}$.

In 2019, Sandu extended MIS methods to higher order accuracy in time, and to
support implicit treatment of the slow components, naming his
newly-proposed methods \emph{Multirate Infinitesimal GARK} (MRI-GARK)
\cite{Sandu2018}.  These allowed the forcing term $r_i$ from \eqref{eq:mis_fast_ode} to be time-dependent and to include the stage solution, $z_i$,
\begin{equation}
  \label{eq:mri_gark_fast_forcing}
  r_i(t) = \frac{1}{\Delta c_i}\sum_{j=1}^{i} \left(\sum_{k=1}^{n_{\Gamma}} \gamma^{\{k\}}_{i,j} \left(\frac{t-T_{i-1}}{\Delta c_i h}\right)^{k-1}\right) \fvec^S(T_{j}, z_j),
\end{equation}
where the coefficients $\Gamma \in \R^{s\times s\times n_{\Gamma}}$ are derived to ensure order conditions for coupling the fast and slow components.

In order to avoid the complication of simultaneously solving for an implicit slow stage $z_i$ while evaluating the fast IVP \eqref{eq:mis_fast_ode}, Sandu required a ``solve-decoupled'' structure, wherein $\gamma^{\{k\}}_{i,i}\ne 0$ only for stages where $\Delta c_i = 0$, allowing the IVP \eqref{eq:mis_fast_ode} to reduce to a standard backward Euler-like implicit solve,
\[
  0 = z_i - z_{i-1} - h \sum_{j=1}^{i} \left(\sum_{k=1}^{n_{\Gamma}} \frac{\gamma_{i,j}^{\{k\}}}{k}\right) \fvec^S(T_{j},z_j).
\]

MRI-GARK methods may similarly include embeddings for temporal adaptivity, by adding an additional row of coupling coefficients such that $\Gamma \in \R^{(s+1)\times s\times n_{\Gamma}}$, and utilizing this row to introduce an embedding stage $(s+1)$ as with MIS methods.  MRI-GARK methods were later extended in \cite{chinomona_implicit-explicit_2021} to support an additional partitioning of $\fvec^S$ into implicit and explicit components,
$\fvec^S = \fvec^I + \fvec^E$, thus enabling implicit-explicit multirate integration for problems of the form \eqref{eq:IVP_3split}.  As with the extension of MIS to MRI-GARK methods, these ImEx-MRI-GARK methods were achieved through modifying the forcing function \eqref{eq:mri_gark_fast_forcing},
\begin{equation}
  \label{eq:imex_mri_gark_fast_forcing}
  \begin{split}
    r_i(t) &= \frac{1}{\Delta c_i}\sum_{j=1}^{i} \left(\sum_{k=1}^{n_{\Gamma}} \gamma^{\{k\}}_{i,j} \left(\frac{t-T_{i-1}}{\Delta c_i h}\right)^{k-1}\right) \fvec^I(T_{j}, z_j) \\
    &+ \frac{1}{\Delta c_i}\sum_{j=1}^{i-1} \left(\sum_{k=1}^{n_{\Omega}} \omega^{\{k\}}_{i,j} \left(\frac{t-T_{i-1}}{\Delta c_i h}\right)^{k-1}\right) \fvec^E(T_{j}, z_j),
  \end{split}
\end{equation}
and modifying the backward Euler-like implicit solve,
\begin{equation}
  \label{eq:imex_mri_gark_implicit_solve}
  0 = z_i - z_{i-1} - h \sum_{j=1}^{i} \left(\sum_{k=1}^{n_{\Gamma}} \frac{\gamma_{i,j}^{\{k\}}}{k}\right) \fvec^I(T_{j},z_j) - h \sum_{j=1}^{i-1} \left(\sum_{k=1}^{n_{\Omega}} \frac{\omega_{i,j}^{\{k\}}}{k}\right) \fvec^E(T_{j},z_j),
\end{equation}
where now the coefficients $\Gamma \in \R^{(s+1)\times s\times n_{\Gamma}}$ and $\Omega \in \R^{(s+1)\times s\times n_{\Omega}}$ are derived to satisfy appropriate order conditions.

We note that ImEx-MRI-GARK methods are a superset of MRI-GARK, which themselves are a superset of MIS methods; thus in the remainder of this paper we will refer to this class of methods as simply MRI.  The accuracy and linear stability analysis for MRI methods (as well as many others) may be understood using the \emph{Generalized-Structure Additive Runge-Kutta (GARK)} framework, developed by Sandu and G{\"u}nther \cite{Sandu2015}.  Since our proposed methods derive from MRI methods, we include the most relevant theoretical results from the above articles here.

\begin{theorem}
  \label{thm:imex_mri_gark_internal_consistency}
  (Thm.~2.2 from \cite{chinomona_implicit-explicit_2021})
  An ImEx-MRI-GARK method fulfills the ``internal consistency'' conditions for any fast method iff
  \[
    \Gamma^{\{1\}}\cdot \mathbf{1} = \Omega^{\{1\}}\cdot \mathbf{1} = \Delta c, \quad\text{and}\quad
    \Gamma^{\{k\}}\cdot \mathbf{1} = \Omega^{\{k\}}\cdot \mathbf{1} = 0 \quad\forall k>1,
  \]
  where $\mathbf{1} \in\R^s$ denotes the vector of all ones, and $\Delta c_i = c_{i}-c_{i-1}$.  If this is satisfied, and both the fast and slow methods have order at least two, then the overall ImEx-MRI-GARK method is second order.
\end{theorem}

While not stated succinctly in these references, we summarize their conclusion for higher order MRI methods as follows.

\begin{theorem}
  \label{thm:mri_method_order}
  An MRI method has order of accuracy $p$ if the solver used for the inner IVP \eqref{eq:mis_fast_ode} has order at least $p$, and the coefficients that define the MRI method satisfy their respective coupling conditions through order $p$.
\end{theorem}

\section{Extended STS methods}
\label{sec:extsts}

We propose a class of ExtSTS methods as MRI methods that replace the fast IVP solve \eqref{eq:mis_fast_ode} with a single step of an arbitrary STS method.  Since STS methods generally have order at most two, then Theorem \ref{thm:mri_method_order} states that the resulting ExtSTS methods will also have order at most two.  Therefore, we focus the following discussion on the development of second-order, embedded, implicit-explicit ExtSTS methods.

Assume we have an ARK method given by implicit and explicit Butcher tables, $\{A^I,b^I,c^I\}$ and $\{A^E,b^E,c^E\}$, where both share the same abscissas, $c^I=c^E = c \in\R^{s}$, with $c_1=0$ and $\Delta c_i = c_{i} - c_{i-1}\ge 0$ for $i=2,\ldots,s$.  Assume further that the method is stiffly accurate, i.e., $c_s=1$, $b^I = A^I_{s,:}$ and $b^E = A^E_{s,:}$, and ``solve-decoupled,'' i.e., $a^I_{i,i} \Delta c_i = 0$, for all $i=2,\ldots,s$.  Then one step from $t_n$ to $t_{n+1} = t_n + h$ of the corresponding ExtSTS method proceeds as
\begin{subequations}
  \label{eq:extsts}
  \begin{align}
    &\notag\text{Let:}\; z_1 := y_n \\
    &\notag \text{For } i = 2,\ldots,s:\\
    &\notag \quad \text{Let:}\: g_i = h\sum\limits_{j=1}^{i-1}\left[ \(a^I_{i,j}-a^I_{i-1,j}\)\, \fvec^R_j + \(a^E_{i,j}-a^E_{i-1,j}\)\, \fvec^A_j\right],\\
    &\notag \quad \text{If } a^I_{i,i} = 0 \text{ and } \Delta c_i > 0:\\
    &\label{eq:extsts_a} \quad
    \begin{cases}
      \text{Let:} &T_{i-1} := t_n + c_{i-1} h, \; T_{i} := t_n + c_{i} h,  \;\text{and}\; v_i(T_{i-1}) := z_{i-1},\\
      \text{Evolve:} &v_i'(t) = \fvec^D \left(t,\, v_i\right) + \frac{g_i}{\Delta c_i h}, \text{ using one STS step over } [T_{i-1}, T_{i}],\; \\
      \text{Let:} &z_{i} := v_i\left(T_{i}\right),
    \end{cases}\\
    &\label{eq:extsts_b} \quad\text{Else, solve for } z_i \text{ from the equation:} \quad z_i - h a^I_{i,i} \fvec^R_i = z_{i-1} + g_i,\\
    &\notag \text{Let}: y_{n+1} := z_{s}.
  \end{align}
\end{subequations}
Here, for brevity we denote $\fvec^R_j := \fvec^R\!\left(T_j,z_j\right)$ and $\fvec^A_j := \fvec^A\!\left(T_j,z_j\right)$.

If embedding coefficients $d^I, d^E \in \mathbb{R}^s$ are provided for the ARK method, then compute the ExtSTS embedding via
\begin{subequations}
  \label{eq:extsts_embedding}
  \begin{align}
    &\notag \text{Let:}\; \tilde{g} = h \sum\limits_{j=1}^{s} \left[ \(d^I_{j}-a^I_{s-1,j}\)\, \fvec^R_j + \(d^E_{j}-a^E_{s-1,j}\)\, \fvec^A_j\right],\\
    &\notag \text{If } \Delta c_s > 0:\\
    &\label{eq:extsts_c}
    \begin{cases}
      \text{Let:} &\tilde{v}(T_{s-1}) := z_{s-1},\\
      \text{Evolve:} &\tilde{v}'(t) = \fvec^D \left(t,\, \tilde{v}\right) + \frac{\tilde{g}}{\Delta c_s h}, \text{ using one STS step over } [T_{s-1},T_{s}],\; \\
      \text{Let:} &\tilde{y}_{n+1} := \tilde{v}\left(T_{s}\right),
    \end{cases}\\
    &\label{eq:extsts_d} \text{Else, set } \tilde{y}_{n+1} = z_{s-1} + \tilde{g}.
  \end{align}
\end{subequations}

\section{Method analysis}
\label{sec:analysis}

Since ExtSTS methods are constructed by combining MRI and STS methods, then Theorem \ref{thm:mri_method_order} ensures that ExtSTS methods may be easily constructed from established MRI methods.

\begin{lemma}
  \label{lem:extsts_order}
  If an ExtSTS method is constructed from an MRI method (explicit, implicit, or ImEx) having order at least $p$, and a STS method (e.g., RKC, RKL, or RKG) of order $p$, then the ExtSTS method will have order $p$.
\end{lemma}

For the specific case above where we built the ExtSTS method directly from an ARK method, we need only verify that the formulas  \eqref{eq:extsts}-\eqref{eq:extsts_embedding} result in a second-order ImEx-MRI-GARK method, in which case Lemma \ref{lem:extsts_order} ensures that the ExtSTS method will also have order two.

\begin{lemma}
  \label{lem:extsts_ark_order}
  Given a pair of Butcher tables $\{A^I,b^I,c^I\}$ and $\{A^E,b^E,c^E\}$ for an ARK method that satisfy the assumptions
  \begin{itemize}
    \item[(a)] each table $\{A^{\sigma},b^{\sigma},c^{\sigma}\}, \sigma=\{I,E\}$ satisfies the second order conditions,\\
      $A^{\sigma}\mathbf{1} = c^{\sigma}$, $b^{\sigma}\cdot\mathbf{1} = 1$, and $b^{\sigma}\cdot c^{\sigma} = \tfrac12$,
    \item[(b)] the tables share abscissa, $c^I = c^E = c$, with $0=c_1 \le c_2 \le \cdots \le c_s = 1$,
    \item[(c)] both tables are stiffly accurate, i.e., $a^{\sigma}_{s,j} = b^{\sigma}_j$ for $\sigma=\{I,E\}$ and $j=1,\ldots,s$,
  \end{itemize}
  then the formulas \eqref{eq:extsts}-\eqref{eq:extsts_embedding} result in a second-order ImEx-MRI-GARK method.
\end{lemma}
\begin{proof}
  Since the formulas \eqref{eq:extsts}-\eqref{eq:extsts_embedding} have the same structure as ImEx-MRI-GARK methods, including the fast IVP time intervals and initial conditions, we need only identify the $\Gamma$ and $\Omega$ coefficients to write the ExtSTS method in ImEx-MRI-GARK form, and verify internal consistency.  The ExtSTS forcing term may be written as
  \begin{align*}
    r_i = \frac{g_i}{\Delta c_i h} = \frac{1}{\Delta c_i} \sum\limits_{j=1}^{i-1}\left[ \(a^I_{i,j}-a^I_{i-1,j}\)\, \fvec^R_j + \(a^E_{i,j}-a^E_{i-1,j}\)\, \fvec^A_j\right].
  \end{align*}
  Identifying the ExtSTS $\fvec^R$ and $\fvec^A$ functions with the ImEx-MRI-GARK $\fvec^I$ and $\fvec^E$ functions, respectively, we have
  \[
    \gamma^{\{1\}}_{i,j} = a^I_{i,j}-a^I_{i-1,j}, \quad \omega^{\{1\}}_{i,j} = a^E_{i,j}-a^E_{i-1,j}, \quad\text{for}\quad i=2,\ldots,s,
  \]
  $\gamma^{\{1\}}_{1,j} = \omega^{\{1\}}_{1,j} = 0$, and $\Gamma^{\{k\}} = \Omega^{\{k\}}=0$ for $k>1$.  Thus for $i=2,\ldots,s$,
  \begin{align*}
    \sum_{j=1}^{s} \gamma^{\{1\}}_{i,j} = \sum_{j=1}^s a^I_{i,j} - \sum_{j=1}^s a^I_{i-1,j} = c^I_i - c^I_{i-1} = \Delta c_i,\\
    \sum_{j=1}^{s} \omega^{\{1\}}_{i,j} = \sum_{j=1}^s a^E_{i,j} - \sum_{j=1}^s a^E_{i-1,j} = c^E_i - c^E_{i-1} = \Delta c_i,
  \end{align*}
  and hence $\Gamma^{\{1\}}\cdot \mathbf{1} = \Omega^{\{1\}}\cdot \mathbf{1} = \Delta c$.  The remaining internal consistency condition, $\Gamma^{\{k\}}\cdot \mathbf{1} = \Omega^{\{k\}}\cdot \mathbf{1} = 0 \quad\forall k>1$ follows trivially from $\Gamma^{\{k\}} = \Omega^{\{k\}}=0$ for $k>1$.

  We further note that the implicit solves when $a^I_{i,i} \ne 0$ and $\Delta c_i = 0$ are also equivalent to those in the ImEx-MRI-GARK method:
  \begin{align*}
    0 \ = \ z_i - h a^I_{i,i} \fvec^R_i - z_{i-1} - g_i
    \ = \ z_i - z_{i-1} - h\sum\limits_{j=1}^{i} \gamma^{\{1\}}_{i,j} \fvec^R_j - h\sum\limits_{j=1}^{i-1} \omega^{\{1\}}_{i,j} \fvec^A_j,
  \end{align*}
  which equals \eqref{eq:imex_mri_gark_implicit_solve}, and where we leveraged the fact that $A^I$ is lower triangular for ARK methods, and thus $a_{i-1,i} = 0$.
\end{proof}

Thus Theorem \ref{thm:mri_method_order} guarantees that one may construct a second order ExtSTS method from any set of second-order ARK tables $\{A^I, b^I, c^I\}$ and $\{A^E, b^E, c^E\}$ that satisfy the assumptions of Lemma \ref{lem:extsts_ark_order}.  We further note that ExtSTS methods may alternately be constructed from any at least second-order MRI method, and that the resulting ExtSTS method will also have order at least two.  Moreover, if STS methods of order $p>2$ are developed in the future, then these may be combined with an order $p$ MRI method to create an ExtSTS method of order $p$.

In the remainder of this section we analyze the linear stability of ExtSTS methods, focusing on those constructed directly from ARK methods using \eqref{eq:extsts}-\eqref{eq:extsts_embedding}.  Similar to \cite{chinomona_implicit-explicit_2021,Sandu2015}, we consider the Dahlquist-like linear test problem
\begin{equation}
  \label{eq:Dahlquist3}
  y'(t) = \lambda^A y + \lambda^D y + \lambda^R y, \quad y(0)=1,
\end{equation}
where we assume that $\lambda^D \in \R^-$ (as is typically assumed for STS methods), while $\lambda^A,\lambda^R \in \mathbb{C}$ are assumed to have negative real part.  We additionally assume that this problem is solved using an ExtSTS method with partitioning $\mathbf{f}^A(t,y) = \lambda^A y$, $\mathbf{f}^D(t,y) = \lambda^D y$, and $\mathbf{f}^R(t,y) = \lambda^R y$.  For the step \eqref{eq:extsts_a} over the interval $[T_{i-1},T_i]$, we consider the second-order Runge--Kutta--Chebyshev (RKC) \cite{sommeijer_rkc_1998} and second-order Runge--Kutta--Legendre (RKL) \cite{meyerSecondorderAccurateSuper2012, meyerStabilizedRungeKutta2014a} methods.
For the purposes of this analysis, we merely assume that these may be represented using a Butcher table $\{A^D,b^D,c^D\}$, such that its number of stages $s_i^D$ is selected as the minimum integer such that the RKC or RKL method will be stable for the IVP
\begin{equation}
  \label{eq:STSDahlquist}
  v_i'(t) = \lambda^D v_i + \frac{g_i}{h_i}, \quad v_i(T_{i-1}) = z_{i-1}
\end{equation}
over a step of size $h_i:=\Delta c_i h$ where $g_i$ is a constant; we note that for RKC this number of stages is \cite{sommeijer_rkc_1998}
\begin{equation}
  \label{eq:RKC_stages}
  s_i^D = \max\(2,\left\lceil\(\dfrac{\Delta c_i h|\lambda^D|}{0.653}\)^{1/2}\right\rceil\)
\end{equation}
and for RKL it is \cite{meyerSecondorderAccurateSuper2012}
\begin{equation}
  \label{eq:RKL_stages}
  s_i^D = \max\(2,\left\lceil\frac{\sqrt{9 + 8 \Delta c_i h|\lambda^D|} - 1}{2}\right\rceil\).
\end{equation}

\begin{lemma}
  \label{lem:STS_stab_function}
  An STS method with Butcher table $\{A^D,b^D,c^D\}$, applied to IVP \eqref{eq:Dahlquist3} for a single step of size $h_i=\Delta c_i h$ to evolve $v_i(T_{i-1}) \to v_i(T_i)$ has the form
  \begin{align*}
    v_i(T_i) &= \left[1 + \Delta c_i \eta^D b^D\cdot (I - \eta_i^D A^D)^{-1}\mathbf{1}\right] v_i(T_{i-1}) \\
    &+ \left[ 1 + \Delta c_i \eta^D b^D\cdot(I - \eta_i^D A^D)^{-1}c^D\right]g_i,
  \end{align*}
  where we denote $\eta^D = h\lambda^D$, and $\mathbf{1}\in\R^{s_i^D}$ is a vector of all ones.  For brevity, we denote this propagation function as
  \begin{equation}
    \label{eq:STS_stab_function}
    v_i(T_i) = \mathcal{R}^D_{RK}(\Delta c_i \eta^D)\, v_i(T_{i-1}) + \mathcal{R}^D_{F}(\Delta c_i \eta^D) g_i.
  \end{equation}
\end{lemma}
\begin{proof}
  Writing the STS method in standard Runge--Kutta form, its internal stages are given by
  \begin{align*}
    &z_j = v_i(T_{i-1}) + h_i \sum_{k=1}^{s_i^D} A^D_{j,k} \(\lambda^D z_k + \frac{g_i}{h_i}\), \quad j=1,\ldots,s_i^D,\\
    \Leftrightarrow&\\
    &z_j - h_i\lambda^D \sum_{k=1}^{s_i^D} A^D_{j,k} z_k = v_i(T_{i-1}) + g_i\sum_{k=1}^{s_i^D} A^D_{j,k}, \quad j=1,\ldots,s_i^D.
  \end{align*}
  Stacking the stages together, $\zvec =
  \begin{bmatrix} z_1 & \cdots & z_{s_i^D}
  \end{bmatrix}^T$, and leveraging the row-sum condition, $A^D_{j,k}\ones = c^D$, we have
  \begin{align*}
    &z_j - \Delta c_i \eta^D \sum_{k=1}^{s_i^D} A^D_{j,k} z_k = v_i(T_{i-1}) + g_i c^D_j, \quad j=1,\ldots,s_i^D\\
    \Leftrightarrow&\\
    &\zvec = \(I - \Delta c_i \eta^D A^D\)^{-1}\ones v_i(T_{i-1}) + \(I - \Delta c_i \eta^D A^D\)^{-1}c^D g_i.
  \end{align*}
  Since the STS method is at least first order then $b^D\cdot\ones = 1$, so the step solution is given by
  \begin{align*}
    v_i(T_i) &= v_i(T_{i-1}) + h_i \sum_{j=1}^{s_i^D} b^D_j \(\lambda^D z_j + g_i\)
    = v_i(T_{i-1}) + \Delta c_i \eta^D b^D\cdot \zvec + g_i,\\
    &= \left[1 + \Delta c_i \eta^D b^D\cdot \(I - \Delta c_i \eta^D A^D\)^{-1}\ones\right] v_i(T_{i-1}) \\
    &\quad + \left[1 + \Delta c_i \eta^D b^D\cdot \(I - \Delta c_i \eta^D A^D\)^{-1}c^D\right] g_i.
  \end{align*}
\end{proof}
We note that when $\Delta c_i=0$, the STS stability function simplifies to
\begin{equation}
  \label{eq:STS_deltaczero_simplification}
  \mathcal{R}^D_{RK}(0)\, v_i(T_{i-1}) + \mathcal{R}^D_F(0)\, g_i = v_i(T_{i-1}) + g_i.
\end{equation}
\begin{lemma}
  \label{lem:extSTS_stab_function}
  The ExtSTS method \eqref{eq:extsts} applied to the problem \eqref{eq:Dahlquist3} has stability function
  \begin{align}
    \label{eq:extsts_stability}
    \mathcal{R}(\eta^D,\eta^A,\eta^R) &= \mathbf{e}_s^T \left[I - P - Q\(\eta^R A^I + \eta^A A^E \)\right]^{-1}\mathbf{e}_1,
  \end{align}
  where $\eta^D = h\lambda^D$, $\eta^A = h\lambda^A$, and $\eta^R = h\lambda^R$, $\mathbf{e}_1$ and $\mathbf{e}_s$ are the first and last columns of the $s\times s$ identity matrix, and the matrices $P$ and $Q$ are given by
  \begin{align}
    \label{eq:Pmatrix}
    P &= \operatorname{diag}\left(
      \begin{bmatrix} \mathcal{R}^D_{RK}(\Delta c_2\eta^D) \\
        \vdots \\ \mathcal{R}^D_{RK}(\Delta c_s\eta^D)
    \end{bmatrix},-1\right)\\
    \label{eq:Qmatrix}
    Q &= \operatorname{diag}\left(
      \begin{bmatrix} 0 \\ \mathcal{R}^D_F(\Delta c_2\eta^D) \\
        \vdots \\ \mathcal{R}^D_F(\Delta c_s\eta^D)
    \end{bmatrix}, 0\right) -
    \operatorname{diag}\left(
      \begin{bmatrix} \mathcal{R}^D_F(\Delta c_2\eta^D) \\
        \vdots \\ \mathcal{R}^D_{F}(\Delta c_s\eta^D)
    \end{bmatrix}, -1\right),
  \end{align}
  where ``diag'' is the standard operator that creates a diagonal matrix from a vector.
\end{lemma}
\begin{proof}
  Leveraging Lemma \ref{lem:STS_stab_function}, the ExtSTS \eqref{eq:extsts} method for the IVP \eqref{eq:Dahlquist3} is
  \begin{align*}
    &\text{Let:}\; z_1 := y_n \\
    &\text{For } i = 2,\ldots,s:\\
    &\quad \text{Let:}\: g_i = \sum\limits_{j=1}^{i-1}\left[ \(A^I_{i,j}-A^I_{i-1,j}\)\, \eta^R + \(A^E_{i,j}-A^E_{i-1,j}\)\, \eta^A\right]z_j,\\
    &\text{Let: } z_{i} =
    \begin{cases}
      \mathcal{R}^D_{RK}(\Delta c_i \eta^D)\, z_{i-1} + \mathcal{R}^D_{F}(\Delta c_i \eta^D)\, g_i, & \text{if } A^I_{i,i} = 0 \text{ and } \Delta c_i > 0,\\
      \left(z_{i-1} + g_i\right)/\left(1 - \eta^R A^I_{i,i}\right), & \text{otherwise}.
    \end{cases}\\
    &\text{Let}: y_{n+1} := z_{s}.
  \end{align*}
  Leveraging \eqref{eq:STS_deltaczero_simplification} and the solve-decoupled condition $A^I_{i,i} \Delta c_i = 0$, this becomes
  \begin{align*}
    &\text{Let:}\; z_1 := y_n \\
    &\text{For } i = 2,\ldots,s:\\
    &\quad \text{Let } g_i = \sum_{j=1}^{i-1} \left[(A^I_{i,j}-A^I_{i-1,j}) \eta^R + (A^E_{i,j}-A^E_{i-1,j}) \eta^A\right] z_j,\\
    &\quad \text{Let:}\; z_i =
    \frac{\mathcal{R}^D_{RK}\left(\Delta c_i \eta^D\right) z_{i-1} + \mathcal{R}^D_{F}(\Delta c_i \eta^D) g_i}{1 - \eta^R A^I_{i,i}},\\
    &\text{Let}: y_{n+1} := z_{s}.
  \end{align*}
  Stacking the stages together into $\zvec =
  \begin{bmatrix} z_1 & \ldots & z_s
  \end{bmatrix}^T$, this becomes
  \begin{align*}
    \zvec &= \mathbf{e}_1 y_n + P \zvec + Q\(\eta^R A^I + \eta^A A^E \)\zvec\\
    y_{n+1} &= \mathbf{e}_s^T \zvec,\\
    \Leftrightarrow\qquad &\\
    y_{n+1} &= \mathbf{e}_s^T \left[I - P - Q\(\eta^R A^I + \eta^A A^E \)\right]^{-1}\mathbf{e}_1 y_n,
  \end{align*}
  Thus the overall ExtSTS propagation function is given by
  \[
    \mathcal{R}(\eta^D,\eta^A,\eta^R) = \mathbf{e}_s^T \left[I - P - Q\(\eta^R A^I + \eta^A A^E \)\right]^{-1}\mathbf{e}_1.
  \]
\end{proof}

Similar to the approach used for ImEx-MRI-GARK methods \cite{chinomona_implicit-explicit_2021}, we examine joint stability regions for a specified ExtSTS method by plotting the region of values $\eta^\sigma$ where $|\mathcal{R}|<1$ under an assumption that the other parameters $\eta^\nu$ for $\nu\ne\sigma$ are taken from an entire region of $\mathbb{C}$.  Specifically, for the full problem \eqref{eq:Dahlquist3} we plot the advective stability regions
\begin{equation}
  \label{eq:extsts_ADR_joint_stability}
  \mathcal{J}_{\theta,\rho} := \left\{\eta^{A} \in \mathbb{C}^{-} \;:\; |\mathcal{R}(\eta^D,\eta^A,\eta^R)|\leq 1,\; \forall \eta^R \in \mathcal{S}_{\theta}^R,\:
  \forall \eta^D \in [-\rho,0)\right\}
\end{equation}
where $\mathcal{S}_{\theta}^{\sigma} := \left\{\eta^{\sigma} \in \mathbb{C}^{-} \;:\; |\arg (\eta^{\sigma})- \pi| \leq \theta\right\}$, and $\eta^D$ is assumed to be real-valued.

We will similarly consider ExtSTS methods for problems where one of $\fvec^R=0$ or $\fvec^A=0$, by plotting the regions
\begin{align}
  \label{eq:extsts_AD_joint_stability}
  \mathcal{J}^A_{\rho} &:= \left\{\eta^{A} \in \mathbb{C}^{-} \;:\; |\mathcal{R}(\eta^D,\eta^A,0)|\leq 1,\; \forall \eta^D \in [-\rho,0) \right\}, \quad\text{and}\\
  \label{eq:extsts_RD_joint_stability}
  \mathcal{J}^R_{\rho} &:= \left\{\eta^{R} \in \mathbb{C}^{-} \;:\; |\mathcal{R}(\eta^D,0,\eta^R)|\leq 1,\; \forall \eta^D \in [-\rho, 0) \right\}.
\end{align}

Finally, for completeness we define the ARK joint stability region for the implicit-explicit splitting of \eqref{eq:Dahlquist3} into $\fvec^I=(\lambda^D+\lambda^R)y=\lambda^Iy$ and $\fvec^E=\lambda^Ay$ as in \cite{zharovsky_class_2015},
\begin{equation}
  \label{eq:ARK_joint_stability}
  \mathcal{J}^{ARK}_{\theta} := \left\{\eta^{E} \in \mathbb{C}^{-} \;:\; |\mathcal{R}^{ARK}(\eta^E,\eta^I)|\leq 1,\; \forall \eta^I \in \mathcal{S}_{\theta}^I \right\},
\end{equation}
where $\mathcal{R}^{ARK}$ is the ARK linear stability function, and $\mathcal{S}_{\theta}^{I}$ is the wedge of angle $2\theta$ in the left complex plane for $\eta^I$, similar to $\mathcal{S}_{\theta}^{\sigma}$ above.

\section{Candidate ExtSTS methods}
\label{sec:methods}

This section constructs a set of proposed second-order ExtSTS methods, and provides analysis of their corresponding stability regions. Following the formulas \eqref{eq:extsts} and \eqref{eq:extsts_embedding}, we construct these methods from embedded ARK Butcher tables, $\{A^I,b^I,c^I,d^I\}$ and $\{A^E,b^E,c^E,d^E\}$.  If either the explicit or implicit Butcher table is zero (i.e., if the underlying method is an explicit Runge--Kutta (ERK) or diagonally-implicit Runge--Kutta (DIRK) method), then the ExtSTS method may also be easily created, under the caveat that the method is applied to a simpler advection-diffusion IVP,
\begin{equation}
  \label{eq:IVP_2split_explicit}
  \yvec'(t) = \fvec^{A}(t,\yvec)+\fvec^{D}(t,\yvec), \quad \yvec(t_0) = \yvec_0,
\end{equation}
or reaction-diffusion IVP,
\begin{equation}
  \label{eq:IVP_2split_implicit}
  \yvec'(t) = \fvec^{R}(t,\yvec)+\fvec^{D}(t,\yvec), \quad \yvec(t_0) = \yvec_0.
\end{equation}

For advection-diffusion-reaction problems, we construct an ImEx ExtSTS method from the embedded ARK method in \cite{giraldo2013implicit}.  Padding to ensure a solve-decoupled and stiffly-accurate structure, our \emph{Giraldo-ARK} method has the Butcher tables

\begin{figure}[htbp]
  \centering
  \label{fig:ARKStability_a}\includegraphics[trim={25 20 40 22}, clip, width=0.48\textwidth]{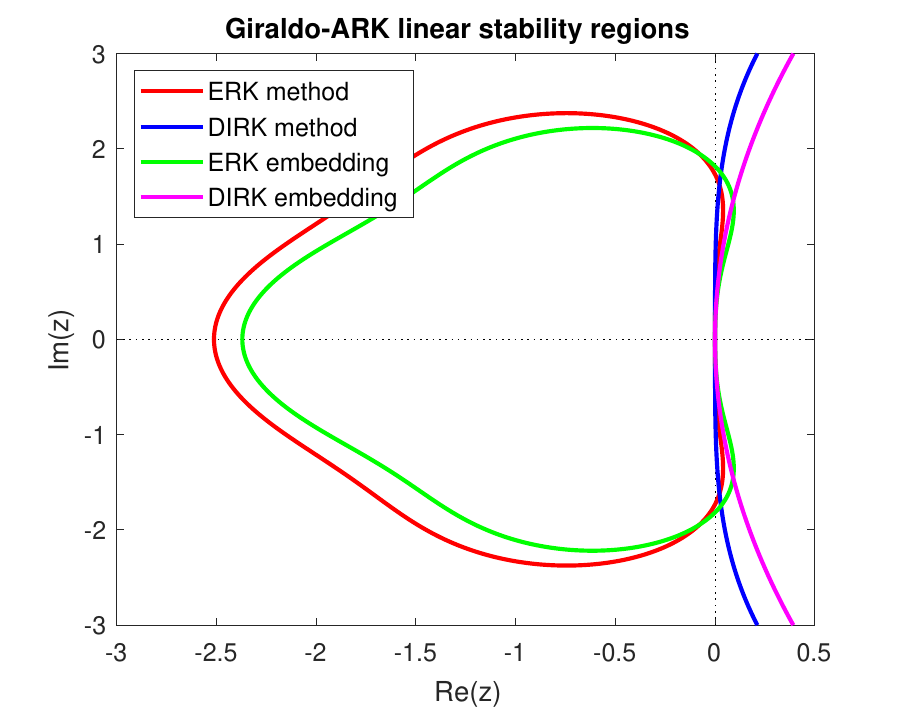}
  \hfill
  \label{fig:ARKStability_b}\includegraphics[trim={40 20 40 22}, clip, width=0.46\textwidth]{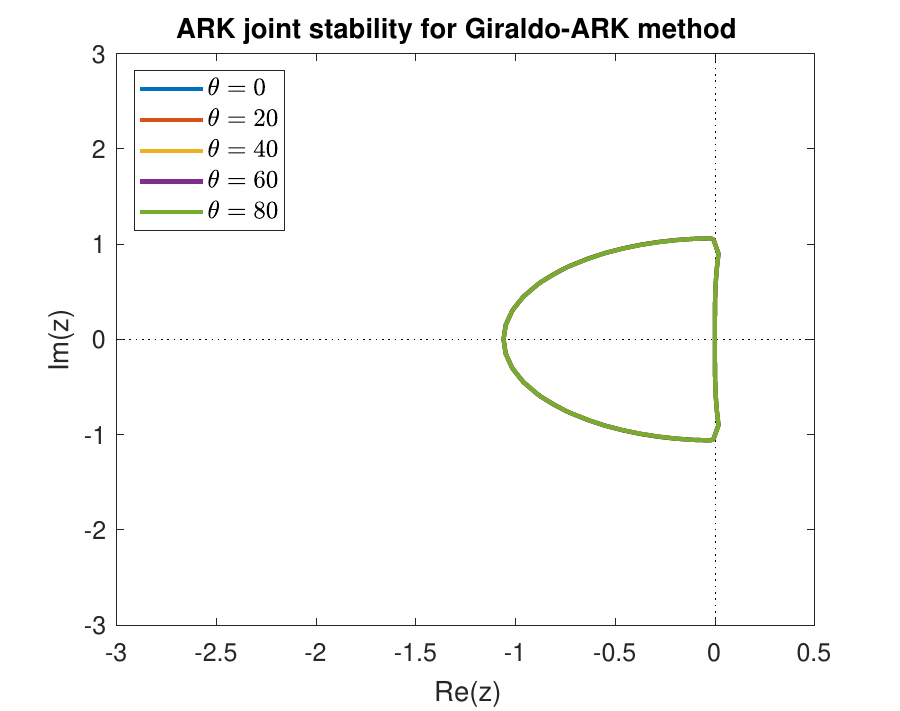}\\
  \vspace{0.5em}
  \label{fig:extsts_stability_a}\includegraphics[trim={25 0 30 22}, clip, width=0.48\textwidth]{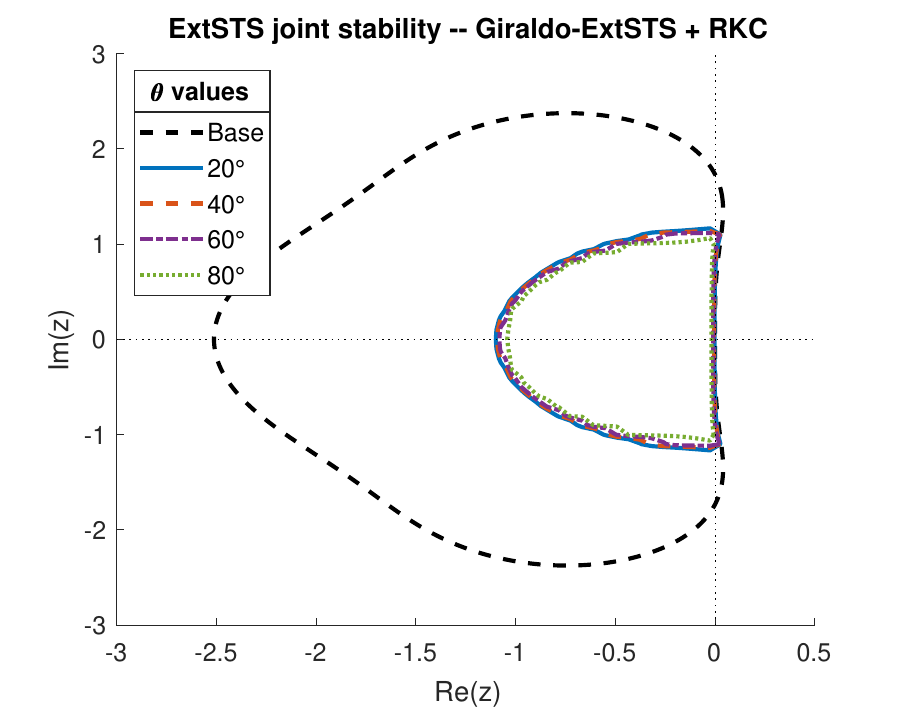}
  \hfill
  \label{fig:extsts_stability_b}\includegraphics[trim={40 0 30 22}, clip, width=0.46\textwidth]{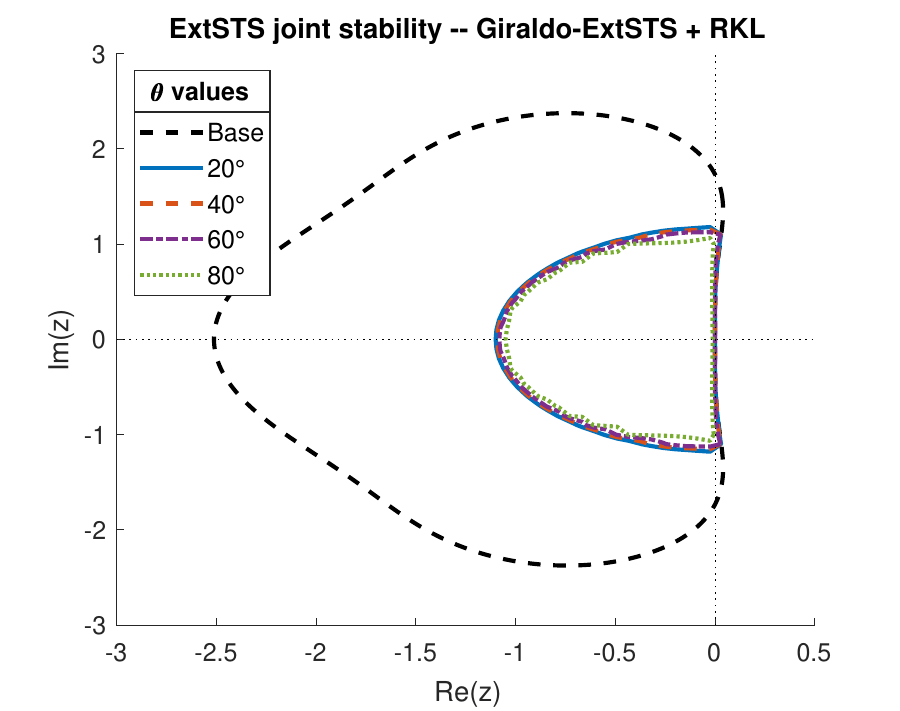}

  \caption{Linear stability regions for the Giraldo-ARK (top) and Giraldo-ExtSTS (bottom) methods.  The top row shows the linear stability of each component table separately (left) and joint ARK stability (right), where we see that the ARK joint stability region $\mathcal{J}^{ARK}_{\theta}$ is slightly smaller than the underlying explicit method, but does not change as the implicit angle $\theta$ varies.  The ExtSTS regions $\mathcal{J}_{\theta,10^6}$ in the bottom row correspond with the RKC (left) and RKL (right) STS sub-methods, showing that the ExtSTS joint linear stability regions remain approximately the same as its ARK joint stability region.}
  \label{fig:ImEx-stability}
\end{figure}

\newcommand{\sqt}{\sqrt{2}}

\begin{align}
  \label{eq:ExtSTS-GiraldoARK}
  &
  \notag
  \renewcommand{\arraystretch}{1.2}
  \begin{array}{c|c|c}
    c & A^E & A^I \\
    \hline
    & d^E & d^I
  \end{array} = \\
  &
  \renewcommand{\arraystretch}{1.2}
  \begin{array}{c|cccccc|cccccc}
    0      & 0                   & 0& 0                   & 0& 0                   & 0
    & 0                   & 0& 0                   & 0& 0                   & 0\\
    2\gamma & 2\gamma             & 0& 0                   & 0& 0                   & 0
    & 2\gamma & 0& 0                   & 0& 0                   & 0\\
    2\gamma & 2\gamma             & 0& 0                   & 0& 0                   & 0
    & \gamma & 0& \gamma & 0& 0                   & 0\\
    1      & \frac{3-2\sqt}{6}   & 0& \frac{3 + 2\sqt}{6} & 0& 0                   & 0
    & \textcolor{red}{0}  & 0& \textcolor{red}{1}  & 0& 0                   & 0\\
    1      & \frac{3-2\sqt}{6}   & 0& \frac{3 + 2\sqt}{6} & 0& 0                   & 0
    & \delta     & 0& \delta     & 0& \gamma & 0\\
    1      & \delta     & 0& \delta     & 0& \gamma & 0
    & \delta     & 0& \delta     & 0& \gamma & 0\\
    \hline
    & \frac{4 - \sqt}{8} & 0& \frac{4 - \sqt}{8} & 0& \delta & 0
    & \frac{4 - \sqt}{8} & 0& \frac{4 - \sqt}{8} & 0& \delta & 0
  \end{array}
\end{align}
where $\gamma = (2-\sqt)/2$ and $\delta = \sqt/4$.  We note that the entries in \textcolor{red}{red} may be modified if desired, so long as internal consistency (i.e., the row-sum condition, $c = A^I\mathbf{1}$) holds.  We note that the implicit portion of this method is A- and L-stable, has a stage order of 1, and has an A-stable embedding.  In \cref{fig:ImEx-stability} we show the linear stability regions for this method, both in isolation (ERK vs DIRK) and its ARK joint stability $\mathcal{J}^{ARK}_{\theta}$, as well as the advective ExtSTS joint stability regions $\mathcal{J}_{\theta,\rho}$ \eqref{eq:extsts_ADR_joint_stability} for $\rho=10^6$ and a range of $\theta$ values for the corresponding \emph{Giraldo-ExtSTS} method when using either RKC or RKL diffusion sub-methods.  Notably, we see that the proposed method shows no reduction in ExtSTS joint stability for either STS method as compared to its ARK joint stability.

\begin{figure}[htbp]
  \centering
  \label{fig:extsts_explicit_stability_a}\includegraphics[trim={25 20 40 22}, clip, width=0.48\textwidth]{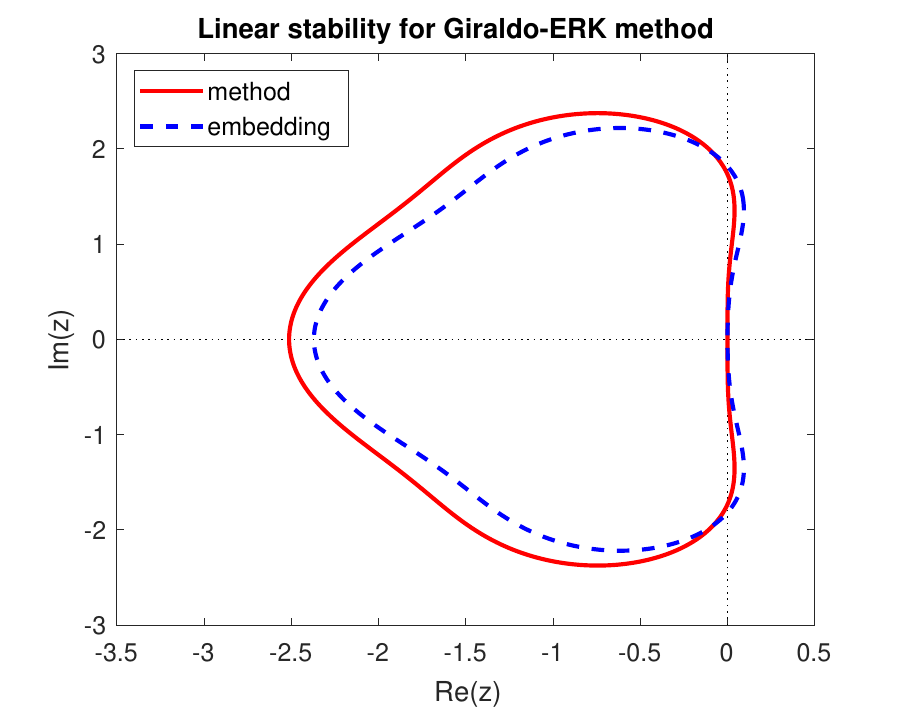}
  \hfill
  \label{fig:extsts_explicit_stability_b}\includegraphics[trim={40 20 40 22}, clip, width=0.46\textwidth]{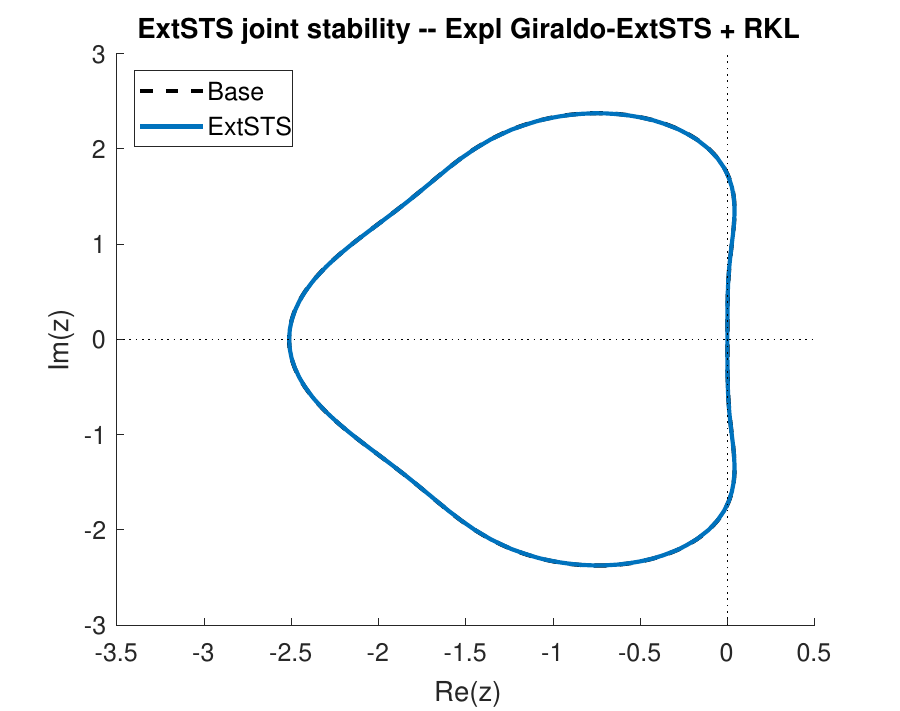}\\
  \vspace{0.5em}
  \label{fig:extsts_explicit_stability_c}\includegraphics[trim={25 0 30 22}, clip, width=0.48\textwidth]{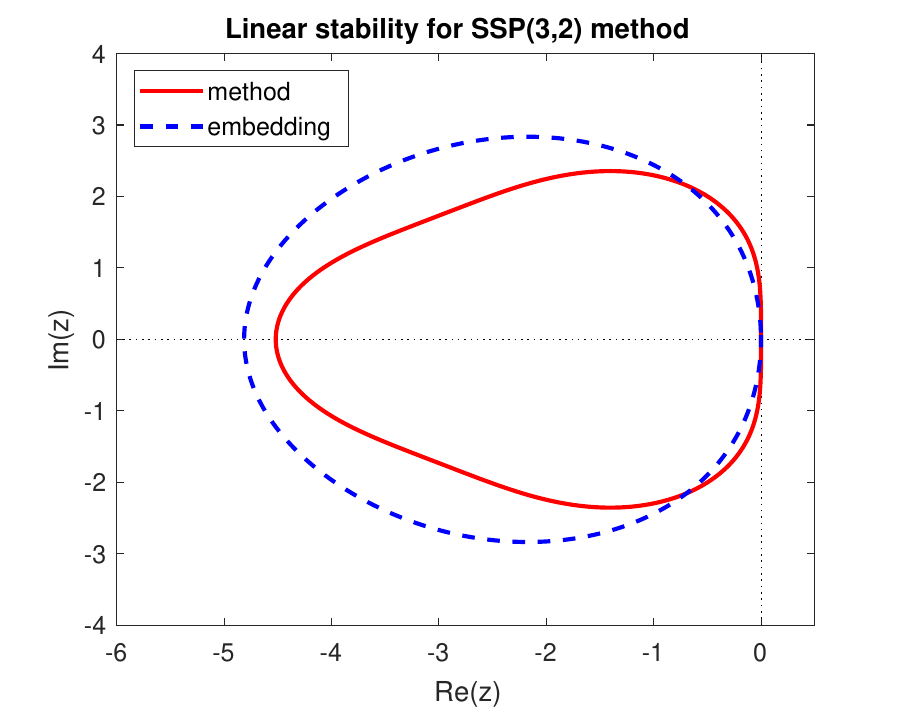}
  \hfill
  \label{fig:extsts_explicit_stability_d}\includegraphics[trim={40 0 30 22}, clip, width=0.46\textwidth]{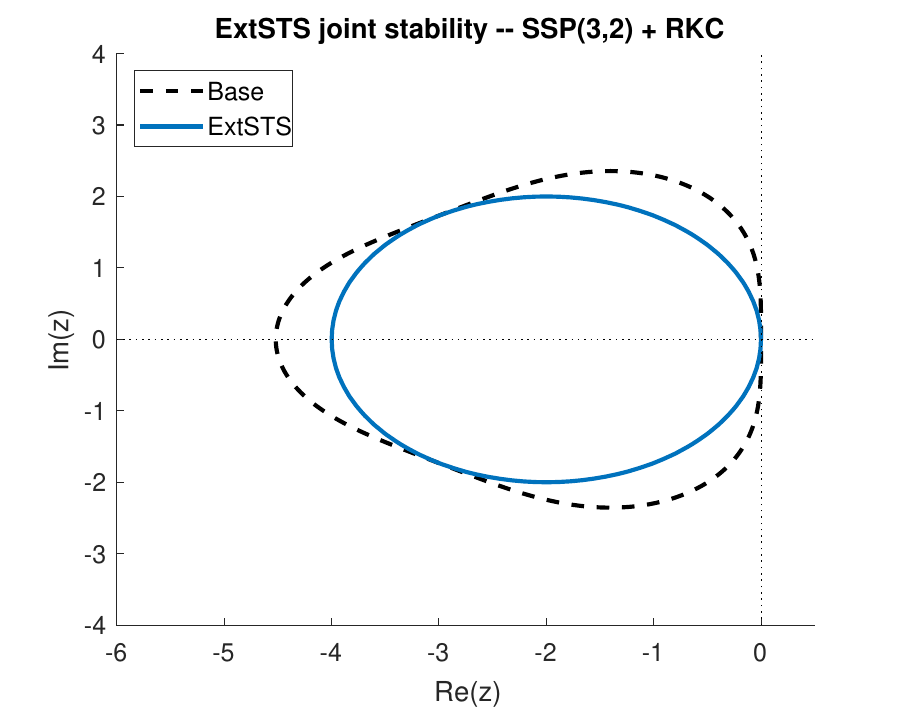}

  \caption{RK Linear stability (left) and ExtSTS joint linear stability $\mathcal{J}^R_{10^6}$ (right) for the explicit Giraldo-ExtSTS method on top, and the SSP32-ExtSTS method on the bottom.}
  \label{fig:extsts_jointstab_explicit}
\end{figure}

For advection-diffusion problems (i.e., where $\fvec^R=0$), we consider explicit ExtSTS methods constructed from two ERK methods, the non-padded explicit table from \eqref{eq:ExtSTS-GiraldoARK}, that we will refer to as \emph{explicit Giraldo-ExtSTS}, and the embedded 3-stage second-order optimal strong-stability-preserving method from \cite{feketeEmbeddedPairsOptimal2022}, \emph{SSP32-ExtSTS}. In stiffly accurate form, these Butcher tables are
\begin{align}
  \label{eq:ExtSTS-GiraldoERK}
  \renewcommand{\arraystretch}{1.2}
  \begin{array}{c|c}
    c & A \\
    \hline
    & d
  \end{array} =
  \renewcommand{\arraystretch}{1.2}
  \begin{array}{c|cccc}
    0 & 0 & 0 & 0 & 0 \\
    2\gamma & 2\gamma & 0 & 0 & 0\\
    1 & \frac{3-2\sqrt{2}}{6} & \frac{3+2\sqrt{2}}{6} & 0 & 0\\
    1 & \delta & \delta & \gamma & 0\\
    \hline
    & \frac{4-\sqrt{2}}{8} & \frac{4-\sqrt{2}}{8} & \delta & 0
  \end{array}
  \qquad\text{and}\qquad
  \renewcommand{\arraystretch}{1.2}
  \begin{array}{c|cccc}
    0 & 0 & 0 & 0 & 0 \\
    \frac12 & \frac12 & 0 & 0 & 0 \\
    1 & \frac12 & \frac12 & 0 & 0 \\
    1 & \frac13 & \frac13 & \frac13 & 0\\
    \hline
    1 & \frac49 & \frac13 & \frac29 & 0
  \end{array}.
\end{align}
As with the ARK method above, we see no difference in linear stability with the ExtSTS methods that use RKC or RKL for the STS sub-solver, so in \cref{fig:extsts_jointstab_explicit} we plot only the original linear stability regions (left) and ExtSTS joint linear stability regions $\mathcal{J}^R_{10^6}$ (right) for these two methods.  Notably, since these do not include an implicit operator, the explicit Giraldo-ExtSTS joint stability region is unchanged from its base Runge--Kutta method, while the SSP32-ExtSTS joint stability region is only slightly smaller than its base counterpart.

\begin{figure}[H]
  \centering
  \label{fig:extsts_implicit_stability_a}\includegraphics[trim={10 0 40 25}, clip, width=0.48\textwidth]{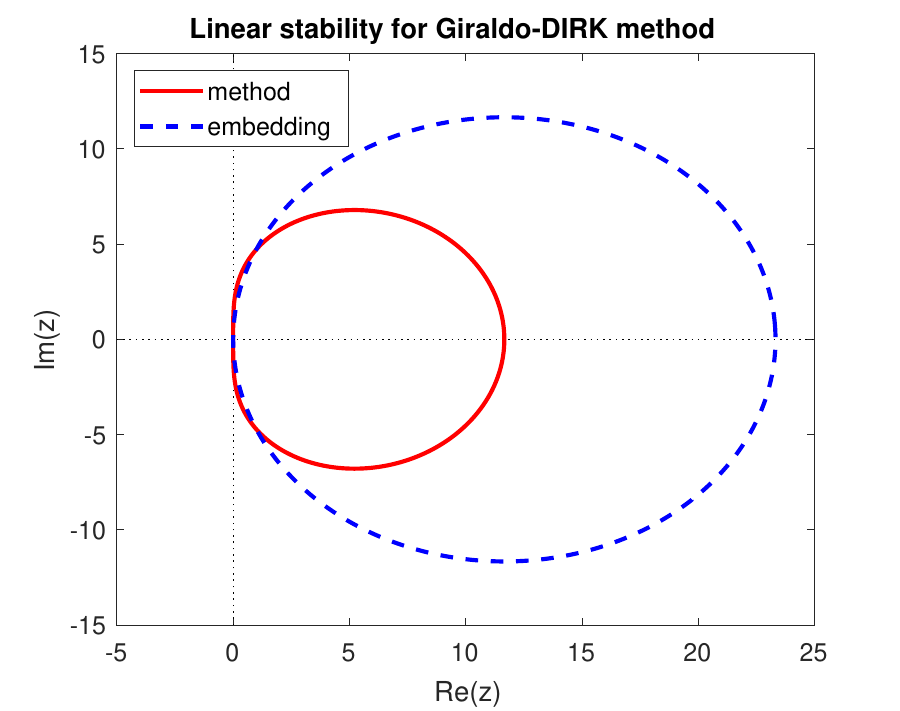}
  \hfill
  \label{fig:extsts_implicit_stability_b}\includegraphics[trim={30 0 40 25}, clip, width=0.46\textwidth]{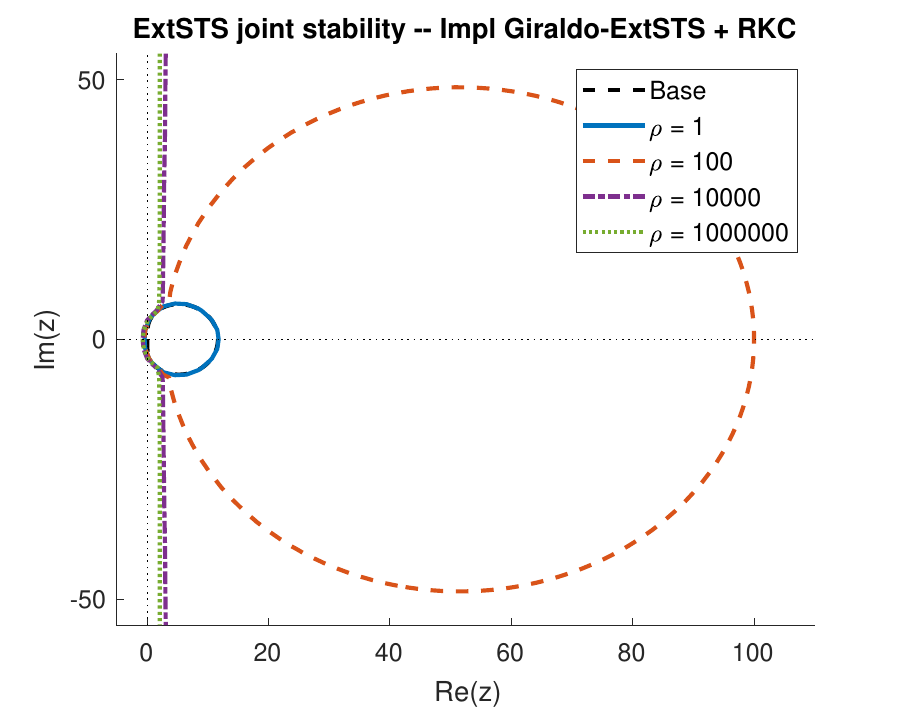}

  \caption{RK Linear stability (left) and ExtSTS joint linear stability $\mathcal{J}^R_{\rho}$ (right) for the implicit Giraldo-ExtSTS method.  We note that $\mathcal{J}^R_{\rho}$ is linearly stable throughout the entire left half-plane for all $\rho$ values.}
  \label{fig:extsts_jointstab_implicit}
\end{figure}

For reaction-diffusion problems (i.e., where $\fvec^A=0$), we consider an implicit ExtSTS method constructed from the DIRK portion of \eqref{eq:ExtSTS-GiraldoARK}, which we will refer to as \emph{implicit Giraldo-ExtSTS}.  Again the choice of STS sub-solver does not affect joint linear stability, so in \cref{fig:extsts_jointstab_implicit} we plot the RK linear stability region (left) and ExtSTS joint linear stability region using RKC (right).  Here, the joint stability region $\mathcal{J}^R_{\rho}$ decreases in size as $\rho$ increases, but the method remains linearly stable throughout the entire left half-plane for all $\rho$ values.

\section{Numerical results}
\label{sec:results}

\begin{figure}[htbp]
  \centering
  \begin{subfigure}[b]{0.49\textwidth}
    \centering
    \includegraphics[trim={0 0 0 20}, clip, width=\textwidth]{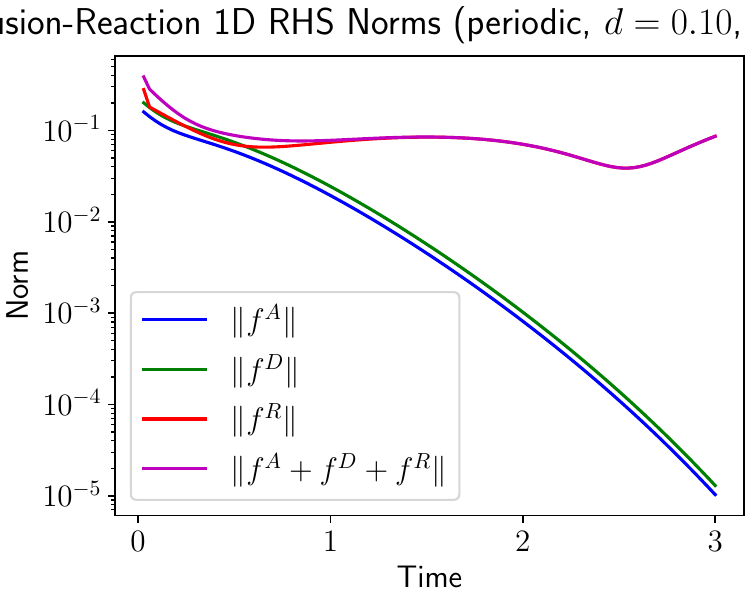}
    \caption{Periodic}
    \label{fig:boundary_effects_a}
  \end{subfigure}
  \hfill
  \begin{subfigure}[b]{0.49\textwidth}
    \centering
    \includegraphics[trim={0 0 0 20}, clip, width=\textwidth]{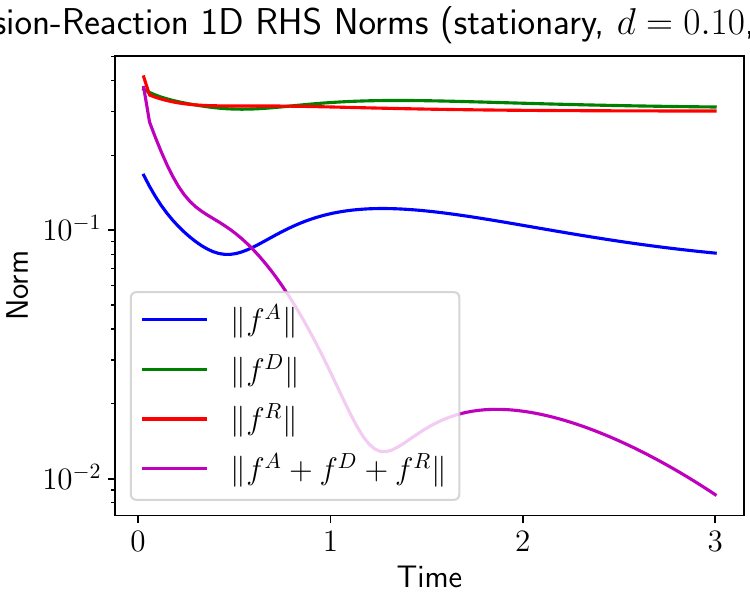}
    \caption{Stationary}
    \label{fig:boundary_effects_b}
  \end{subfigure}
  \caption{Norms of the advection, diffusion, and reaction operators for the stiff Brusselator problem \eqref{eq:adr_1D} with periodic (left) and stationary (right) boundary conditions, for the parameters $c=0.5$, $d=0.1$, $r=1$, and $\epsilon=10^{-2}$.  In the periodic case, there is very little cancellation between the operators, while in the stationary case there is significant cancellation, particularly as the problem progresses in time.}
  \label{fig:boundary_effects}
\end{figure}

We examine the performance of the proposed methods on both one- and two-dimensional systems of advection-diffusion-reaction equations, based on the well-known ``Brusselator'' benchmark problem by Prigogine \cite{Prigogine1967}.  While this is a simple benchmark problem, it captures the essential features of more complex multiphysics systems that couple advection, diffusion, and reaction processes.  The one-dimensional problem is given by:
\begin{equation}
  \label{eq:adr_1D}
  \begin{split}
    \partial_t u &= -c \partial_x u + d \partial_{xx} u + r\(A - (w + 1) u + v u^2\),\\
    \partial_t v &= -c \partial_x v + d \partial_{xx} v + r\(w u - v u^2\),\\
    \partial_t w &= -c \partial_x w + d \partial_{xx} w + r\(\tfrac{B - w}{\epsilon}- w u\),
  \end{split}
\end{equation}
over the domain $(t,x) \in [0,3]\times[0,1]$, discretized on a regular spatial grid with 512 nodes using second-order finite-differences, where $u$, $v$, and $w$ represent the concentrations of chemical species over this domain, $c$ is the advection speed, $d$ is the diffusion strength, $\epsilon$ is the stiffness parameter, and the species with constant concentration over time are $A = 1$ and $B = 3$.  The baseline values of these coefficients are $c = 0.5$, $d = \{0.01, 0.1\}$, $r=1$, and $\epsilon=\{10^{-2}, 10^{-4}, 10^{-6}\}$, although these may be changed in specific tests that follow.  We split the operators in these equations in their natural components, with $\mathbf{f}^A$ containing the advection terms $-c\partial _x (\cdot)$, $\mathbf{f}^D$ containing the diffusion terms $d\partial_{xx}(\cdot)$, and $\mathbf{f}^R$ containing the reaction terms grouped with the $r$ parameter.  The initial conditions are given by
\begin{equation}
  \label{eq:adr_1D_initial}
  \begin{split}
    u(0,x) &=  A  + 0.1 \sin(2 \pi x),\\
    v(0,x) &= \tfrac{B}{A} + 0.1 \sin(2 \pi x),\\
    w(0,x) &=  B  + 0.1 \sin(2 \pi x).
  \end{split}
\end{equation}
We consider two versions of this problem.  The first uses periodic boundary conditions and the second uses stationary boundary conditions, i.e., $u_t(t,0) = u_t(t,1) = v_t(t,0) = v_t(t,1) = w_t(t,0) = w_t(t,1) = 0$.  Of these, we consider the second to be more challenging -- since the boundary values are held fixed the operators $\fvec^A$, $\fvec^D$ and $\fvec^R$ to oppose one another near the boundary, stressing the coupling between these operators in the time integration method.  This may be observed in Figure \ref{fig:boundary_effects}, where we separately plot each of $\|\fvec^A\|$, $\|\fvec^D\|$, $\|\fvec^R\|$, and $\|\fvec^A + \fvec^D + \fvec^R\|$ as functions of time, showing that in the periodic case there is very little cancellation between the operators, whereas in the stationary case there is significant cancellation.

In the following tests, we compare the proposed ExtSTS methods against a set of competing multi-physics time integration methods that are standard for problems with this structure.  The most direct comparison is against the PIROCK solver.  To perform those tests, we modified a test problem that was provided in the original Fortran-77 PIROCK code to mirror the physical problems above; for each test we run PIROCK with the same configuration (implicit/explicit partitioning of the problem, tolerances, etc.) as we use for the ExtSTS methods.  We also compare against standard ImEx additive Runge--Kutta methods that are formed using the non-padded versions of the Butcher table \eqref{eq:ExtSTS-GiraldoARK},
\begin{align}
  \label{eq:GiraldoARK}
  \renewcommand{\arraystretch}{1.5}
  \begin{array}{c|c|c} c & A^E & A^I\\ \hline & d^E & d^I
  \end{array}
  \; = \;
  \renewcommand{\arraystretch}{1.5}
  \begin{array}{c|cccc|cccc}
    0      & 0                   & 0                   & 0       & 0
    & 0                   & 0                   & 0       & 0\\
    2\gamma & 2\gamma            & 0                   & 0       & 0
    & \gamma              & \gamma              & 0       & 0\\
    1      & \frac{3-2\sqt}{6}   & \frac{3 + 2\sqt}{6} & 0       & 0
    & \delta     & \delta     & \gamma  & 0\\
    1      & \delta     & \delta     & \gamma  & 0
    & \delta     & \delta     & \gamma  & 0\\
    \hline
    & \frac{4 - \sqt}{8} & \frac{4 - \sqt}{8} & \delta & 0
    & \frac{4 - \sqt}{8} & \frac{4 - \sqt}{8} & \delta & 0
  \end{array},
\end{align}
where $\gamma$ and $\delta$ are the same as in \eqref{eq:ExtSTS-GiraldoARK}.  For this ARK method, $\mathbf{f}^A$ is treated explicitly, while $\mathbf{f}^D+\mathbf{f}^R$ is treated implicitly.  When performing fixed step-size tests, we also compare a second order Strang--Marchuk operator splitting method \cite{Marchuk1968, Strang1968}, wherein $\mathbf{f}^D$ is treated using either the RKC or RKL method, and the remaining $\mathbf{f}^A$ and $\mathbf{f}^R$ terms are treated using an additive Runge--Kutta method with the table \eqref{eq:GiraldoARK} above.  In each case, all implicit solves (reactions for ExtSTS, PIROCK, and Strang, and reaction+diffusion for ARK) are performed using a Newton method with a direct linear solver.  All computational tests in this paper are provided in an accompanying GitHub repository \cite{CEDADemonstrationsRepo}.

\subsection{Super-time-stepping sub-solver tests}
\label{sec:sts_tests}

Before comparing solver performance between algorithm families, we first examine the sensitivity of the proposed ExtSTS methods and the Strang splitting method to the choice of STS sub-solver applied to the $\fvec^D$ operator.  We perform tests using temporal adaptivity with fixed absolute tolerance $\text{atol}=10^{-11}$ and relative tolerances $\text{rtol}=\{10^{-2.5}, 10^{-3}, \ldots, 10^{-6.5}\}$, as well tests with fixed step sizes $h=\{1/40, 1/80, \ldots, 1/5120\}$.  We compare all results against a reference solution obtained using the fourth order ARK4(3)7L[2]SA$_1$ adaptive method from \cite{KenCarp:19}, and with a tighter relative tolerance of $10^{-10}$.  We separately consider three configurations of the problem \eqref{eq:adr_1D}: advection-diffusion-reaction, advection-diffusion, and reaction-diffusion, with stationary boundary conditions.  For the advection-diffusion and reaction-diffusion problems, we also include results using ExtSTS methods based on the \emph{ERK22a} and \emph{ESDIRK34a} MRI-GARK methods from \cite{Sandu2018}.

\begin{figure}
  \centering
  \includegraphics[trim={0 20 330 20}, clip, width=0.274\textwidth]{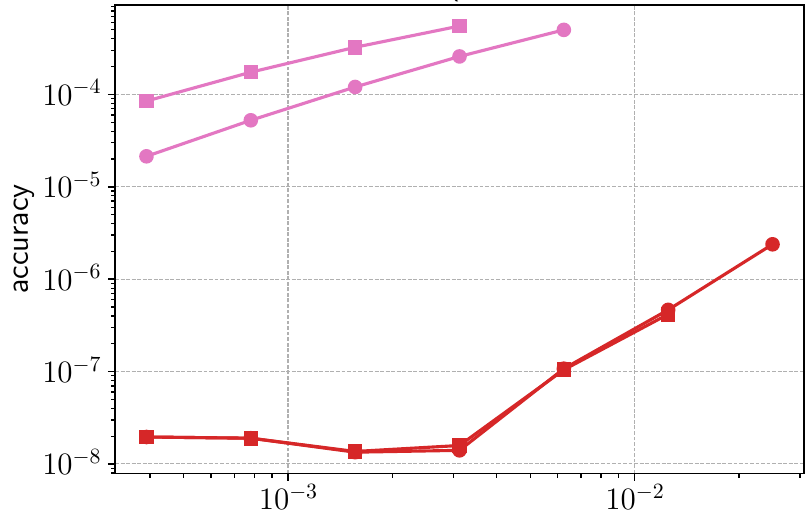}
  \includegraphics[trim={20 20 330 20}, clip, width=0.26\textwidth]{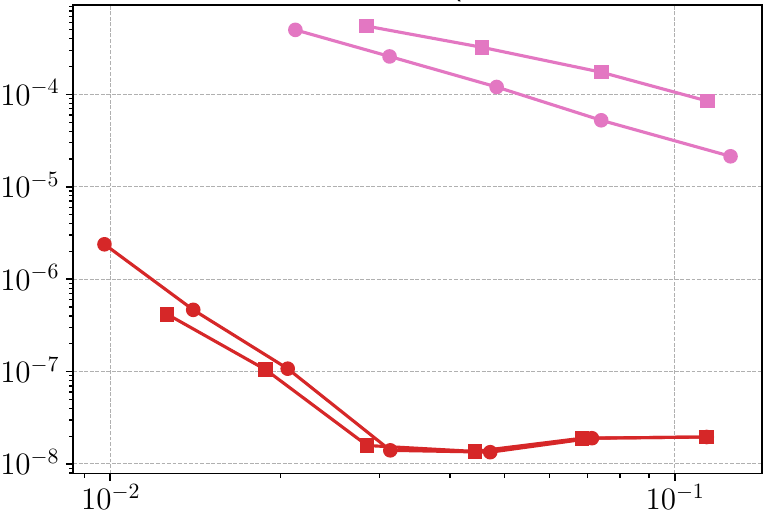}
  \includegraphics[trim={20 20 330 20}, clip, width=0.26\textwidth]{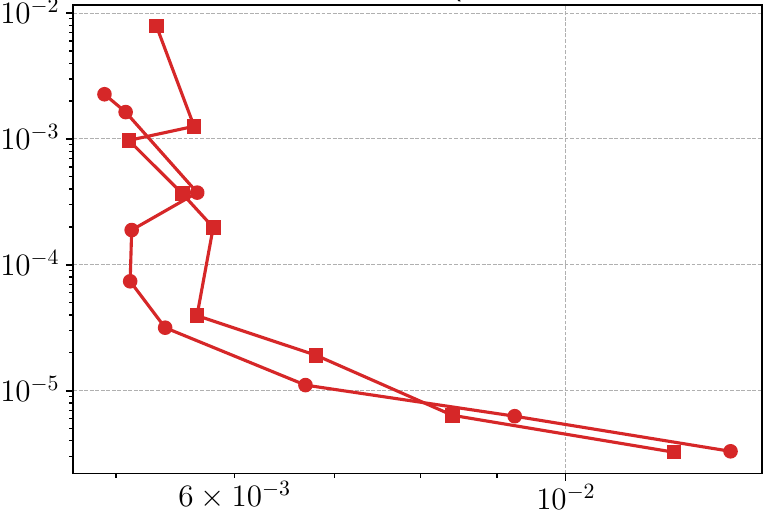}
  {\includegraphics[trim={400 50 100 20}, clip, width=0.18\textwidth]{adr1D_stationary_fixed_runtime_efficiency_RKLvRKC.pdf}}\\
  \includegraphics[trim={0 20 330 20}, clip, width=0.274\textwidth]{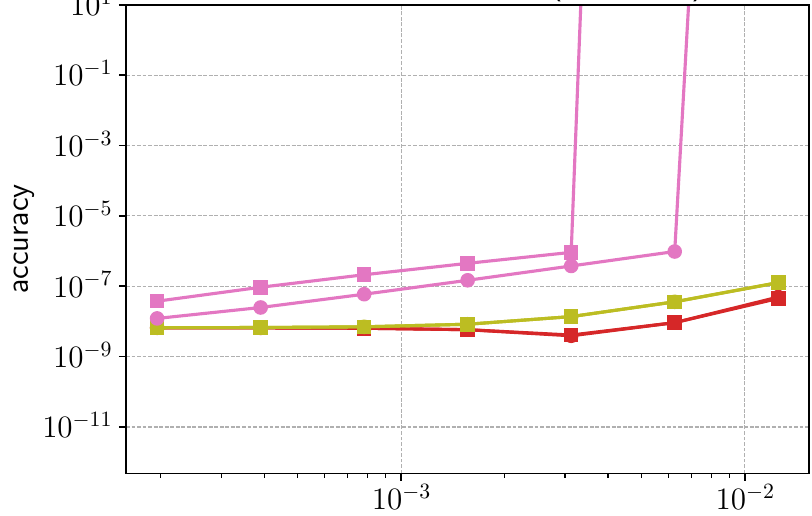}
  \includegraphics[trim={20 20 330 20}, clip, width=0.26\textwidth]{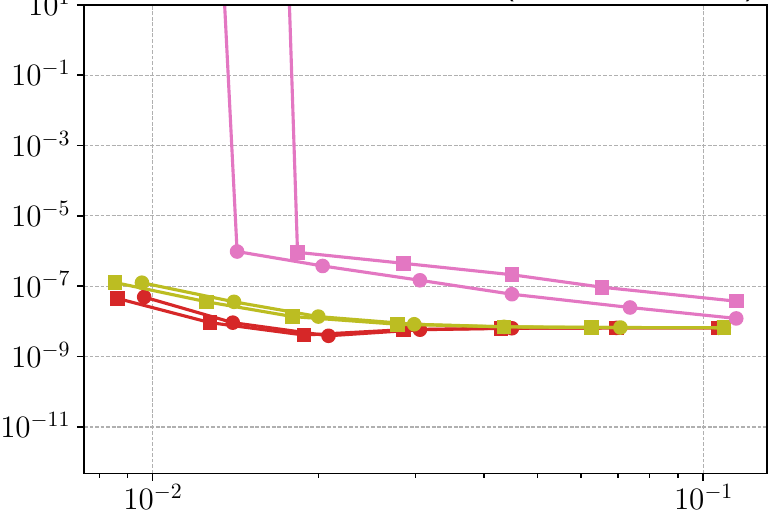}
  \includegraphics[trim={20 20 330 20}, clip, width=0.26\textwidth]{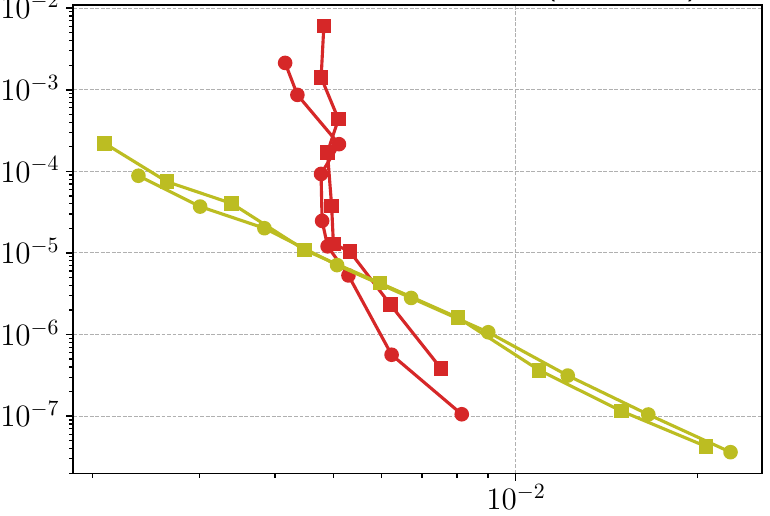}
  {\includegraphics[trim={400 50 100 20}, clip, width=0.18\textwidth]{ad1D_stationary_fixed_runtime_efficiency_RKLvRKC.pdf}}\\
  \includegraphics[trim={0 0 330 20}, clip, width=0.274\textwidth]{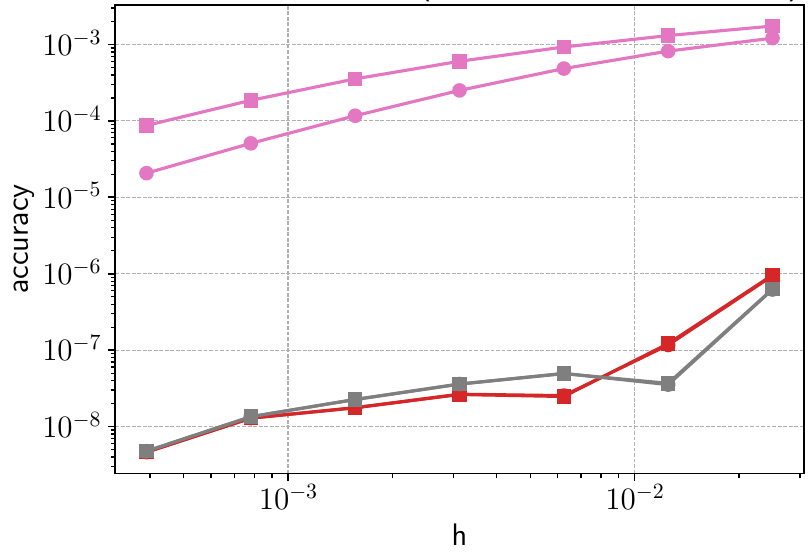}
  \includegraphics[trim={20 0 330 20}, clip, width=0.26\textwidth]{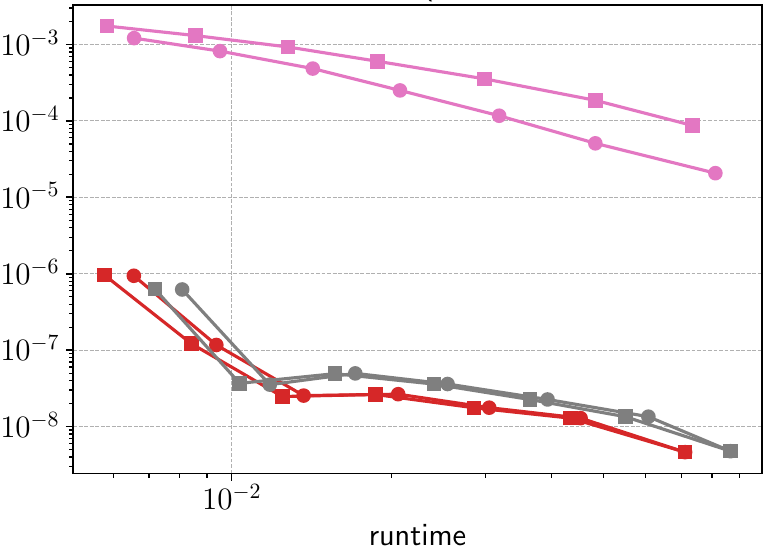}
  \includegraphics[trim={20 0 330 20}, clip, width=0.26\textwidth]{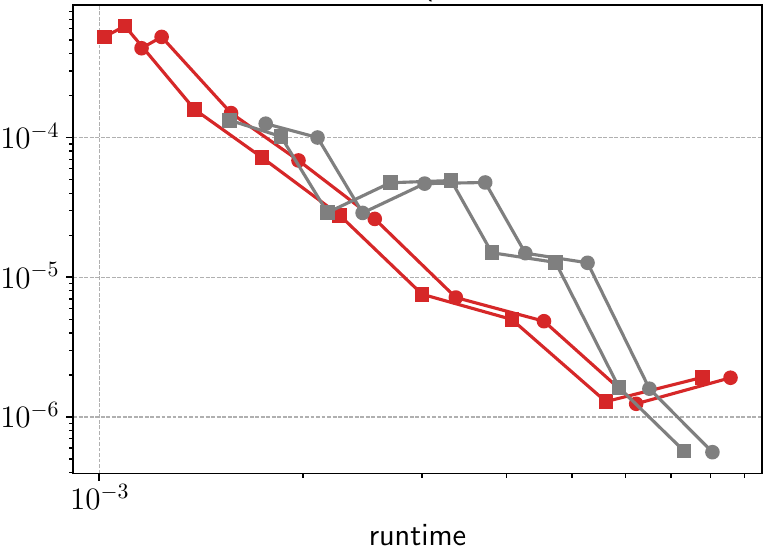}
  {\includegraphics[trim={400 50 100 20}, clip, width=0.18\textwidth]{rd1D_stationary_fixed_runtime_efficiency_RKLvRKC.pdf}}\\
  \caption{Convergence and runtime efficiency of the ExtSTS and Strang methods using RKC and RKL for the diffusion sub-solver.  Methods using RKC are plotted using square markers, while methods using RKL are plotted using circular markers.  Top row: the advection-diffusion-reaction problem \eqref{eq:adr_1D} with stationary boundary conditions, $d=0.1$, $c=0.5$, $r=1$, and $\epsilon=10^{-2}$.  Middle row: the advection-diffusion problem \eqref{eq:adr_1D} with stationary boundary conditions, $d=0.1$, $c=0.5$, and $r=0$.  Bottom row: the reaction-diffusion problem \eqref{eq:adr_1D} with stationary boundary conditions, $d=0.1$, $c=0.0$, and $\epsilon=10^{-2}$. The fixed-step tests use step sizes, $h=\{1/40, 1/80, \ldots, 1/2560\}$.  The first and last rows use fixed step sizes, $h=\{1/40, 1/80, \ldots, 1/2560\},$ and the second row uses $h=\{1/80, 1/80, \ldots, 1/5120\}$.}
  \label{fig:sts_tests}
\end{figure}

In the Figure \ref{fig:sts_tests} we plot the fixed-step convergence, fixed-step runtime efficiency, and adaptive-step runtime efficiency for these tests.  By ``runtime efficiency,'' we mean the ability of a method to achieve a certain level of accuracy in the least amount of computational time. We assess this by plotting the solution accuracy as a function of computational time, to the end that the most efficient method in a given comparison has values closer to the bottom-left of the plot.
We see that in each case, ExtSTS method convergence is unaffected by the choice of RKC or RKL, whereas the Strang solver strongly prefers RKL over RKC.  Due to RKC's slightly larger stability region than RKL, the ExtSTS methods using RKL require slightly more computational effort per step, leading to a slight edge in runtime efficiency for RKC when using fixed time steps; however, this benefit largely disappears when using adaptive time steps.  Thus in all subsequent tests that compare the performance of ExtSTS methods against both Strang and PIROCK, we use RKL for the diffusion sub-solver in the ExtSTS and Strang methods.  We note that the PIROCK solver always uses the ROCK2 STS method for its diffusion component.

\subsection{Advection-diffusion-reaction tests in 1D}
\label{sec:adr1d_tests}

\begin{figure}
  \centering
  \includegraphics[trim={0 0 330 20}, clip, width=0.274\textwidth]{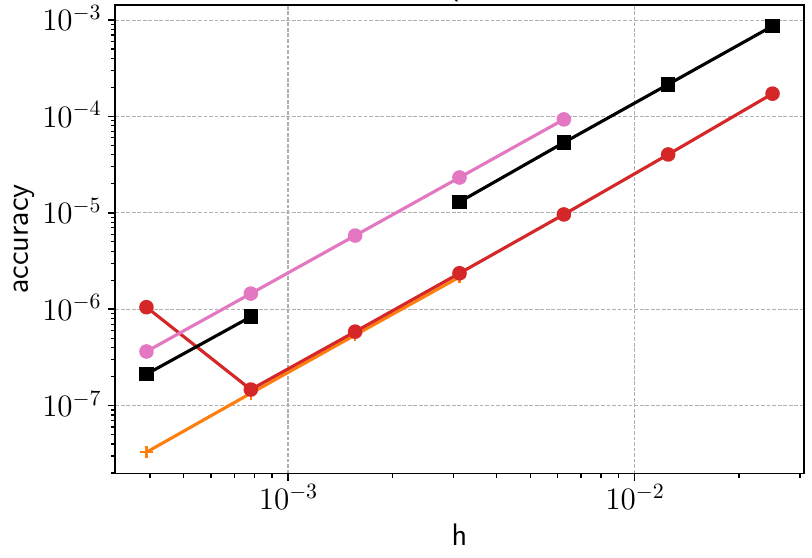}
  \includegraphics[trim={20 0 330 20}, clip, width=0.26\textwidth]{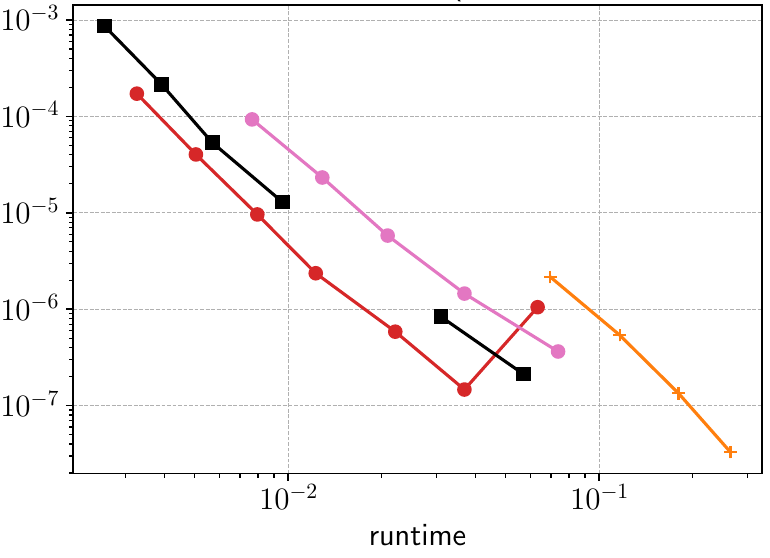}
  \includegraphics[trim={20 0 330 20}, clip, width=0.26\textwidth]{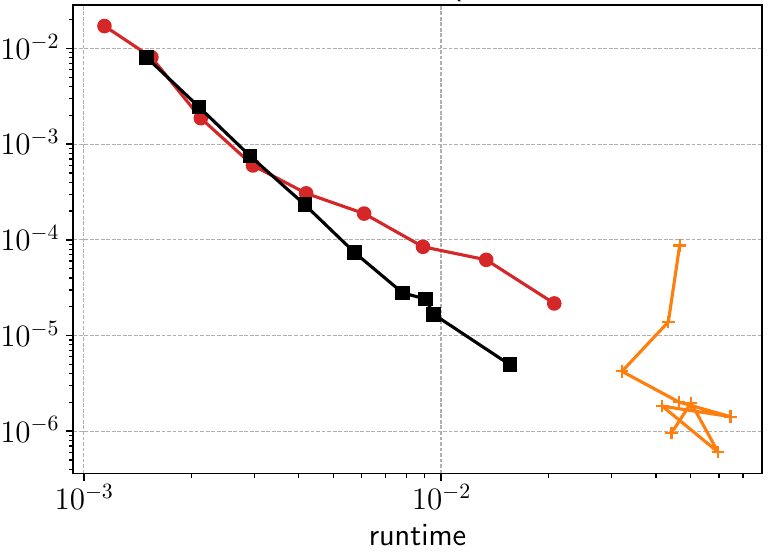}
  \includegraphics[trim={400 60 150 20}, clip, width=0.15\textwidth]{adr1D_periodic_fixed_runtime_efficiency_d0.01_eps1.0e-04.pdf}\\
  \includegraphics[trim={0 0 330 20}, clip, width=0.274\textwidth]{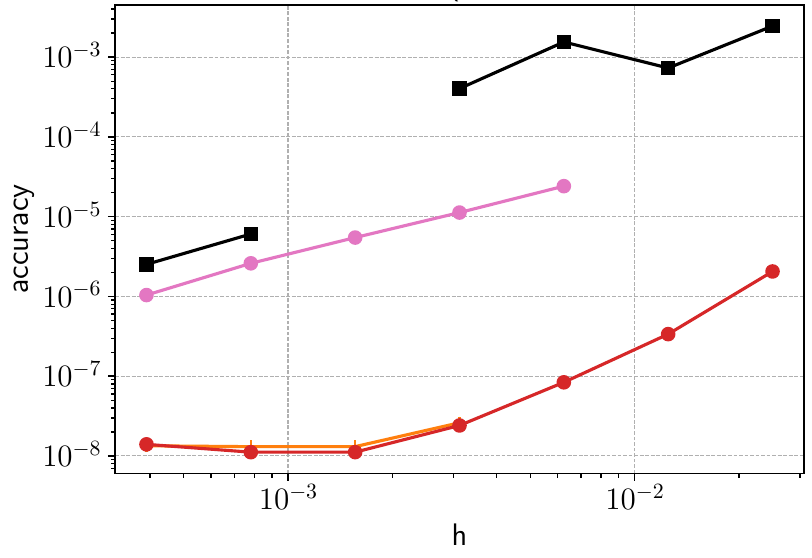}
  \includegraphics[trim={20 0 330 20}, clip, width=0.26\textwidth]{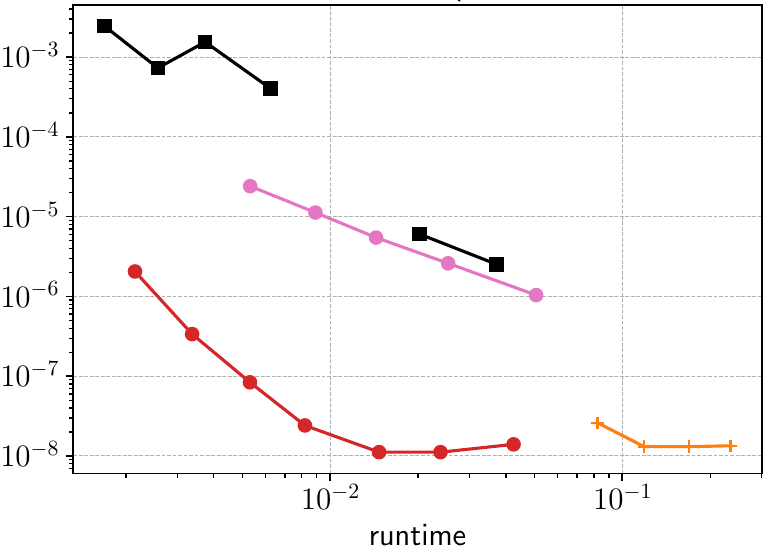}
  \includegraphics[trim={20 0 330 20}, clip, width=0.26\textwidth]{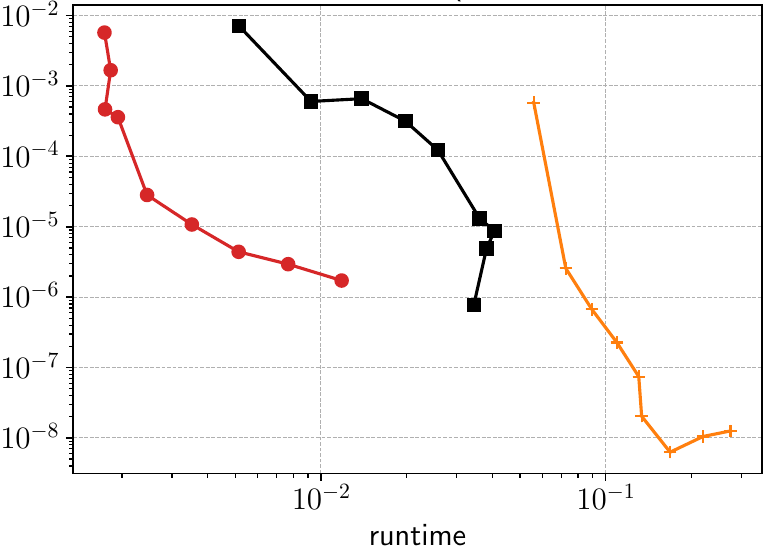}
  \includegraphics[trim={400 0 150 120}, clip, width=0.15\textwidth]{adr1D_stationary_fixed_runtime_efficiency_d0.01_eps1.0e-04.pdf}\\
  \includegraphics[trim={0 0 330 20}, clip, width=0.274\textwidth]{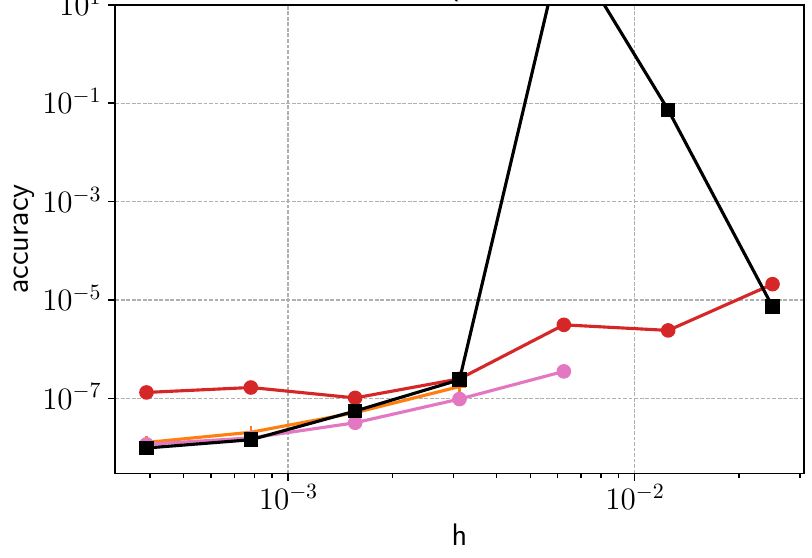}
  \includegraphics[trim={20 0 330 20}, clip, width=0.26\textwidth]{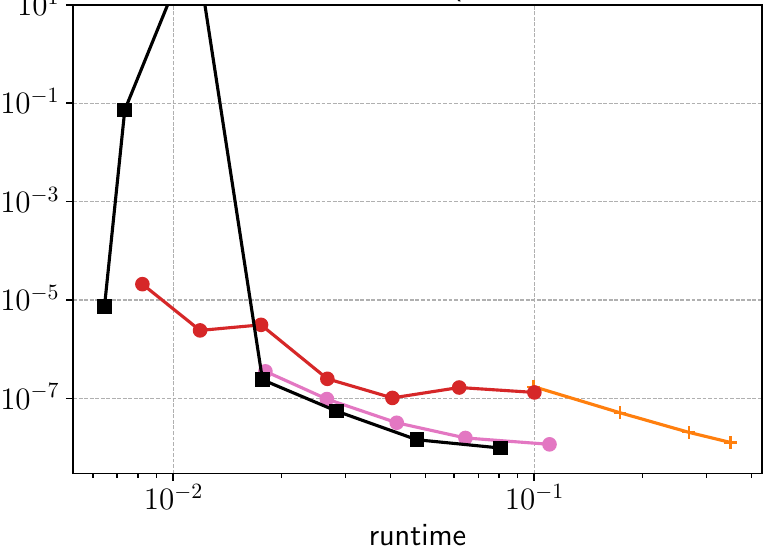}
  \includegraphics[trim={20 0 330 20}, clip, width=0.26\textwidth]{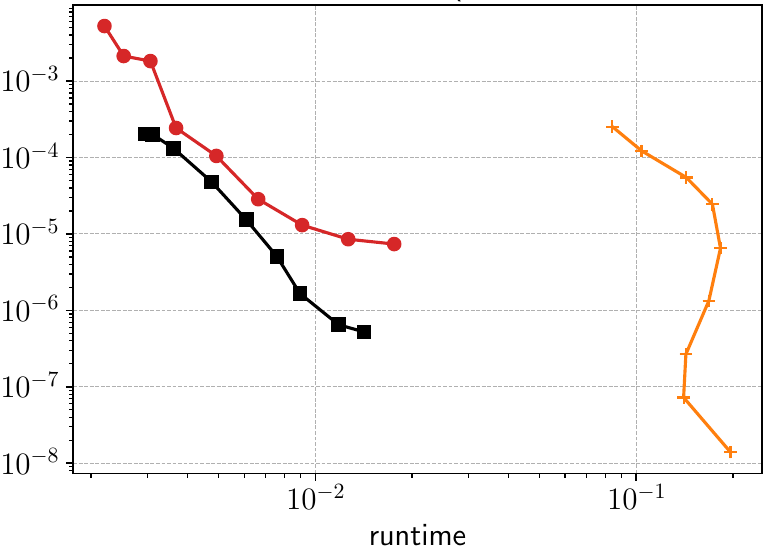}
  \includegraphics[trim={400 0 150 120}, clip, width=0.15\textwidth]{adr1D_periodic_fixed_runtime_efficiency_d0.10_eps1.0e-02.pdf}\\
  \includegraphics[trim={0 0 330 20}, clip, width=0.274\textwidth]{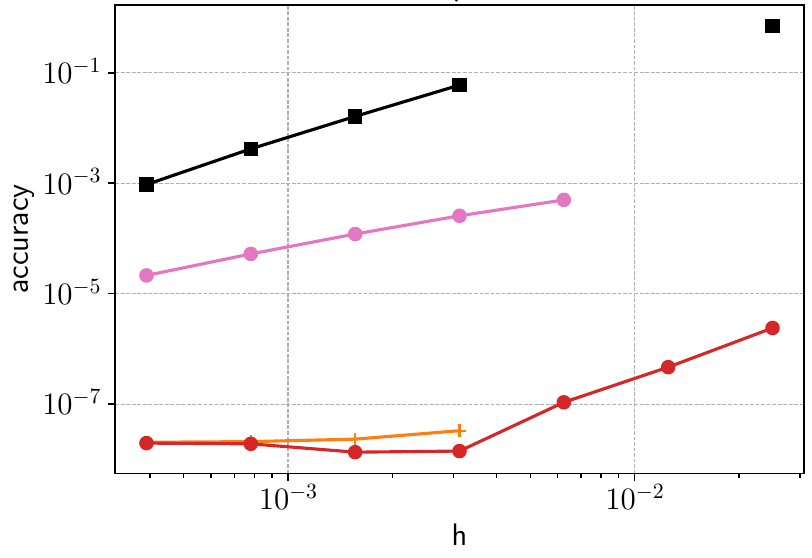}
  \includegraphics[trim={20 0 330 20}, clip, width=0.26\textwidth]{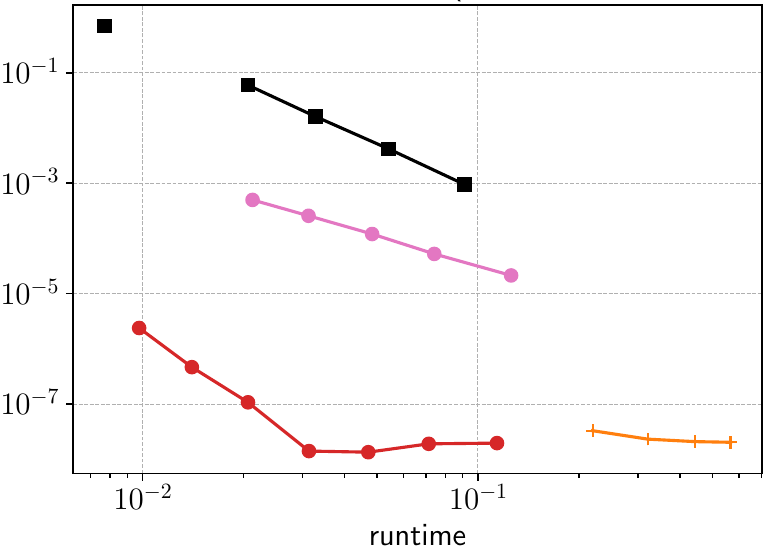}
  \includegraphics[trim={20 0 330 20}, clip, width=0.26\textwidth]{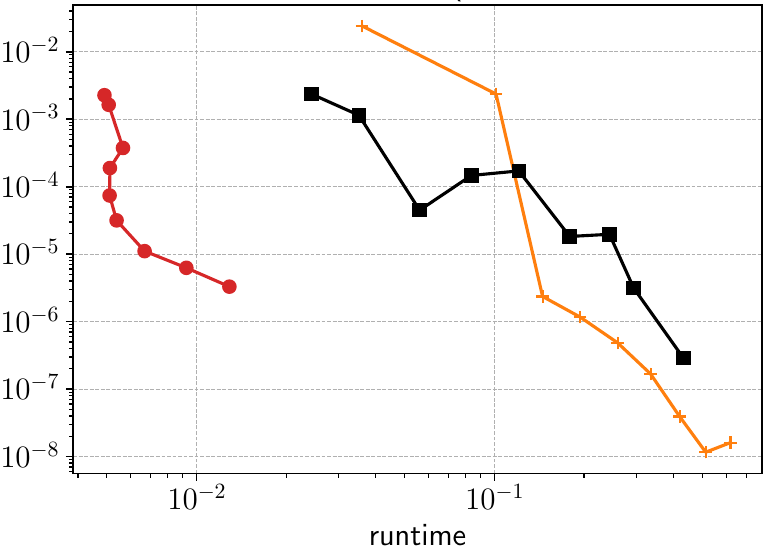}
  \includegraphics[trim={400 0 150 120}, clip, width=0.15\textwidth]{adr1D_stationary_fixed_runtime_efficiency_d0.10_eps1.0e-02.pdf}
  \caption{Convergence (left), runtime efficiency of the fixed (middle) and adaptive (right) methods on the 1D advection-diffusion-reaction test problem with parameters $d=0.01$ and $\varepsilon=10^{-4}$ (rows 1 and 2), and with the parameters $d=0.1$ and $\varepsilon=10^{-2}$ (rows 3 and 4).  Periodic boundary conditions are shown in rows 1 and 3, and stationary boundary conditions in rows 2 and 4.}
  \label{fig:adr1d}
\end{figure}

Our first comparisons focus on the full advection-diffusion-reaction problem \eqref{eq:adr_1D}.  We perform tests both using temporal adaptivity and fixed step sizes, using the same tolerances and step sizes as in Section \ref{sec:sts_tests}.  In Figure \ref{fig:adr1d} we present the convergence and runtime efficiency of the adaptive and fixed step methods on this problem with both periodic and stationary boundary conditions, for the pairs of reaction and diffusion parameters $(\varepsilon,d)=\{(10^{-4},0.01), (10^{-2}, 0.10)\}$, corresponding to problems dominated by reaction and diffusion, respectively.

Due to the tight coupling between each of $\mathbf{f}^A$, $\mathbf{f}^D$, and $\mathbf{f}^R$ within ARK methods, it is unsurprising that they are among the most accurate as fixed step sizes $h$ are reduced.  Since each of ExtSTS, PIROCK, and Strang have the effect of evaluating $f^A$ and $f^R$ less frequently than $f^D$, their accuracy is generally slightly worse than ARK; however their looser coupling generally results in significantly better runtime efficiency than ARK.

Also unsurprisingly, when using fixed time steps, it is clear that the selection of step sizes is critical to performance.  Notably, for both sets of boundary conditions and sets of parameters, only the ExtSTS methods are able to robustly compute solutions at all requested step sizes; both ARK and Strang are generally unstable at the largest step sizes, and PIROCK fails for some intermediate step sizes.  We additionally note that both PIROCK and Strang have significantly larger error for the problems with stationary boundary conditions than periodic, whereas ARK and ExtSTS are largely insensitive to the boundary condition type.  We hypothesize that this is due to the weak coupling between the advection and diffusion components within each PIROCK step, wherein the advection terms are only performed after all diffusion terms have been computed (as is done in operator splitting).  As previously mentioned, unlike in the periodic case, the stationary boundary conditions create a scenario where the advection and diffusion terms are in direct opposition at the downwind boundary, thus methods which do not sufficiently couple these operators are at a disadvantage.

The inconsistent accuracy for the Strang and PIROCK methods, as opposed to the consistent ExtSTS method accuracy, causes significant differences in both fixed step and adaptive runtime efficiency performance.  Here, for periodic boundary conditions PIROCK is generally slightly more efficient than the other methods, but it cannot compete for the more challenging stationary boundary conditions, where ExtSTS methods are by far the most efficient.

\subsection{Advection-diffusion tests in 1D}
\label{sec:ad1d_tests}

To predict the performance of these methods on advection-diffusion problems without implicit reactions, we use our baseline coefficients but disable reactions (i.e., $r=0$).  We perform the same tests as above but using the explicit Giraldo-ExtSTS and SSP32-ExtSTS methods \eqref{eq:ExtSTS-GiraldoERK}.  For this set of tests, only the ARK methods require implicit solves, and they again use a Newton method with a direct linear solver.

\begin{figure}
  \centering
  \includegraphics[trim={0 20 330 20}, clip, width=0.38\textwidth]{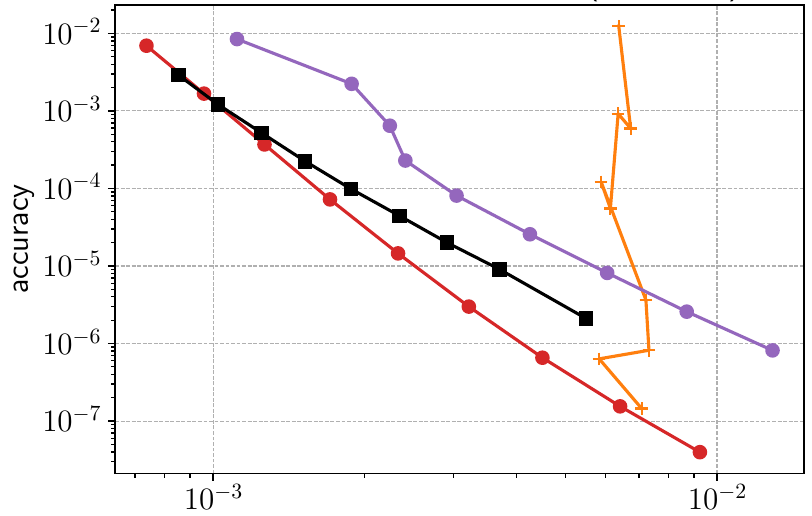}
  \includegraphics[trim={20 20 330 20}, clip, width=0.36\textwidth]{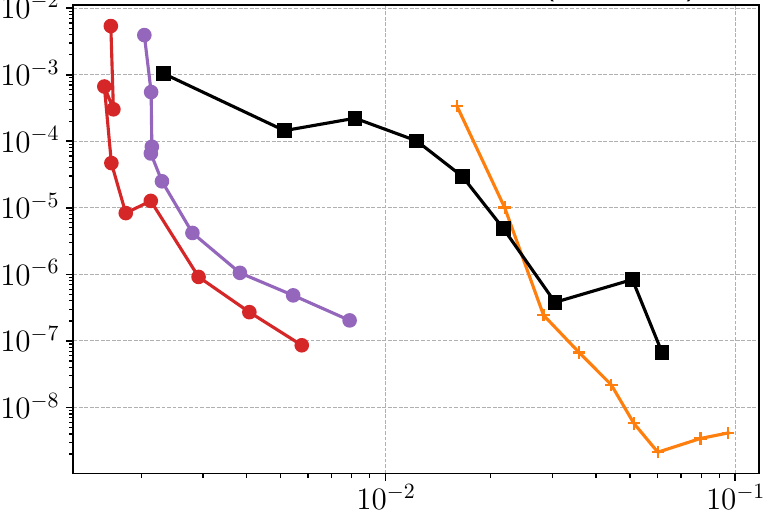}
  \includegraphics[trim={400 30 150 20}, clip, width=0.17\textwidth]{ad1D_stationary_adaptive_runtime_efficiency_d0.01.pdf}\\
  \includegraphics[trim={0 0 330 20}, clip, width=0.38\textwidth]{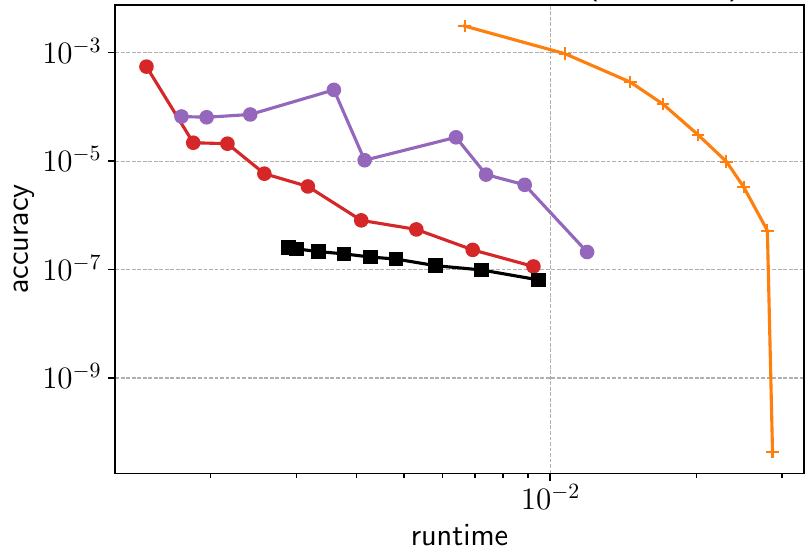}
  \includegraphics[trim={20 0 330 20}, clip, width=0.36\textwidth]{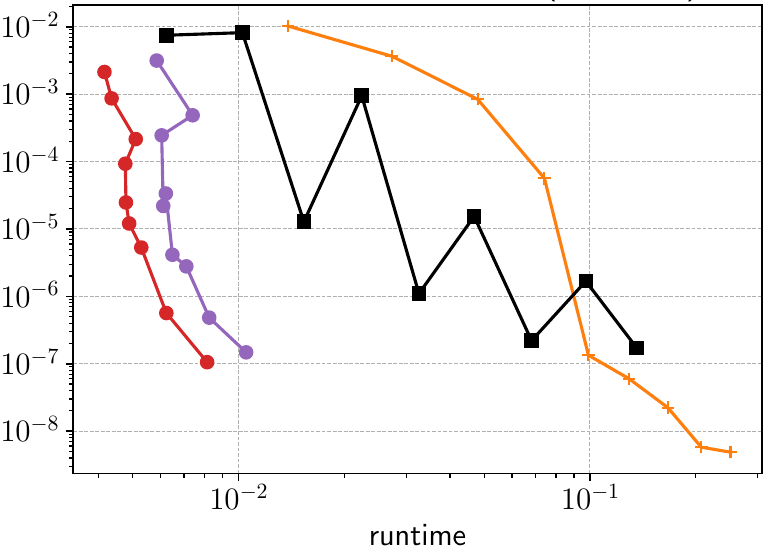}
  \includegraphics[trim={400 0 150 120}, clip, width=0.17\textwidth]{ad1D_stationary_adaptive_runtime_efficiency_d0.01.pdf}\\
  \caption{Adaptive runtime efficiency of the 1D advection-diffusion test problem with diffusion parameters $d=0.01$ (top) and $d=0.10$ (bottom), using periodic boundary conditions (left) and stationary boundary conditions (right).}
  \label{fig:ad1d}
\end{figure}

In Figure \ref{fig:ad1d} we show the runtime efficiency of the adaptive methods on the advection-diffusion test problem with both periodic and stationary boundary conditions for the diffusion parameters $d=\{0.01, 0.10\}$.  We note that with reactions removed, nearly all fixed-step methods exhibit their expected convergence rates, so those plots are omitted.

The methods compare similarly to the advection-diffusion-reaction results in Figure \ref{fig:adr1d}, in that nearly all of the adaptive methods are able to achieve a wide range of solution accuracy.  Performance-wise, the ExtSTS methods perform consistently for both boundary condition types, whereas PIROCK methods only perform well for the simpler periodic boundary condition case.  We note that in these tests, the ExtSTS method formed using the Giraldo-ERK method is uniformly more efficient than the ExtSTS method formed using the SSP32-ERK method; however, the difference is generally small, so for applications with shocks that may require SSP methods for the advection components, we anticipate that SSP32-ExtSTS may be a suitable choice.  Furthermore, since the ExtSTS and PIROCK methods do not require any implicit solves for this problem, they are significantly more efficient than the ARK methods, with the performance gap widening as the diffusion parameter $d$ is increased.  While we do not show the fixed-step results here, and thus the Strang method is not shown, we note that its runtime performance relative to both Strang and ExtSTS is largely similar to the advection-diffusion-reaction problem Figure \ref{fig:adr1d}, in that it performs somewhat better than PIROCK for stationary boundary conditions and somewhat worse than PIROCK for periodic boundary conditions, but in all cases was worse than ExtSTS methods.

\subsection{Reaction-diffusion tests in 1D}
\label{sec:rd1d_tests}

\begin{figure}
  \centering
  \includegraphics[trim={0 20 330 20}, clip, width=0.38\textwidth]{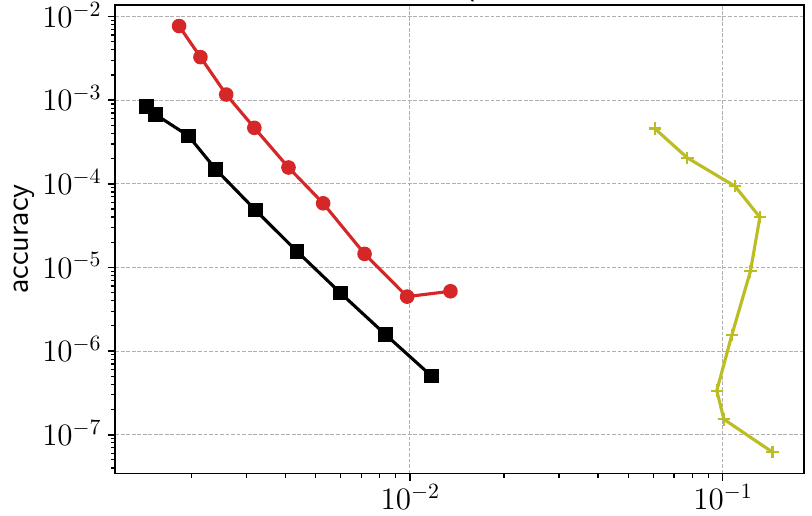}
  \includegraphics[trim={20 20 330 20}, clip, width=0.36\textwidth]{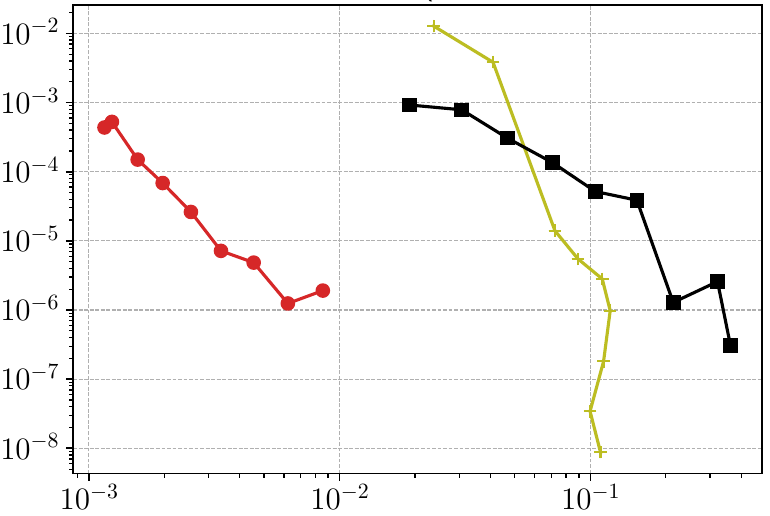}
  \includegraphics[trim={400 30 150 20}, clip, width=0.17\textwidth]{rd1D_stationary_adaptive_runtime_efficiency_d0.10_eps1.0e-02.pdf}\\
  \includegraphics[trim={0 0 330 20}, clip, width=0.38\textwidth]{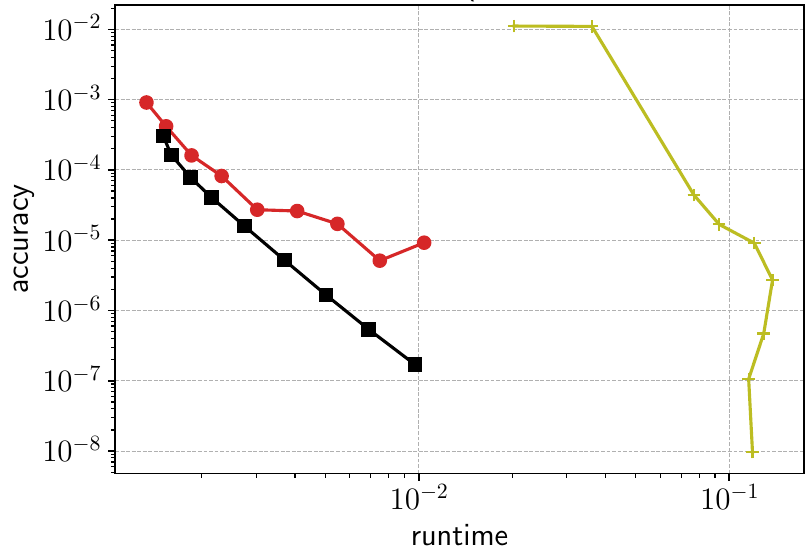}
  \includegraphics[trim={20 0 330 20}, clip, width=0.36\textwidth]{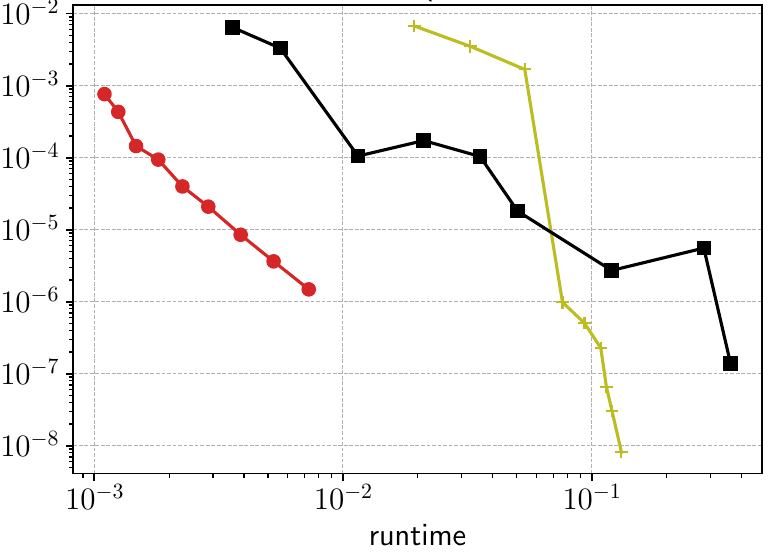}
  \includegraphics[trim={400 0 150 120}, clip, width=0.17\textwidth]{rd1D_stationary_adaptive_runtime_efficiency_d0.10_eps1.0e-04.pdf}
  \caption{Runtime efficiency of the 1D reaction-diffusion test problem with parameters $d=0.10$ and $\varepsilon=10^{-2}$ (top row) and with parameters $d=0.10$ and $\varepsilon=10^{-4}$ (bottom row), using periodic (left) and stationary (right) boundary conditions.}
  \label{fig:rd1d}
\end{figure}

Our final one-dimensional tests examine performance of our proposed methods on reaction-diffusion problems, so we run using our baseline coefficients but now disable advection (i.e., $c=0$).  We compare the performance of the implicit Giraldo-ExtSTS method against a standard diagonally-implicit Runge--Kutta method using the implicit portion of \eqref{eq:ExtSTS-GiraldoARK}, and the PIROCK method.  Here, all methods require implicit solves (reactions for ExtSTS and PIROCK, or reaction+diffusion for ARK), where we again use a Newton method with a direct linear solver.  Finally, since reaction-diffusion problems are frequently quite stiff, we consider the parameters $\varepsilon=\{10^{-2}, 10^{-4}\}$.  We again focus our result on runtime efficiency as the parameters are varied, juxtaposing the cases with periodic and stationary boundary conditions in Figure \ref{fig:rd1d}.

Here, the PIROCK methods are the most efficient option for the configuration using periodic boundary conditions, although the ExtSTS method is a close second place.  These roles again switch for the configuration using stationary boundary conditions, since the ExtSTS method performance remains unchanged but PIROCK performance deteriorates significantly.  Additionally, we note that the ARK methods are generally much less efficient for this problem, arising from the fact that, they are unable to achieve the same level of accuracy for the more challenging stationary boundary conditions.  The ExtSTS methods are able to achieve a wide range of solution inaccuracies for both boundary condition types, and are significantly more efficient than the ARK methods in all cases.

\subsection{Advection-diffusion-reaction tests in 2D}
\label{sec:adr2d_tests}

To assess whether the previous one dimensional results carry over to two dimensions, we conduct a series of tests on a 2D advection-diffusion-reaction problem from the original PIROCK paper \cite{abdullePIROCKSwissknifePartitioned2013a}:
\begin{equation}
  \label{eq:adr_2D}
  \begin{split}
    \partial_t u &= d \nabla^2 u + \mathbf{c}_u\cdot\nabla u + A + u^2 v - (B + 1) u,\\
    \partial_t v &= d \nabla^2 v + \mathbf{c}_v\cdot\nabla v + B u - u^2 v,
  \end{split}
\end{equation}
over the domain $(t,x)\in[0, 5]\times[0, 1]^2$, discretized on a $400\times 400$ regular spatial grid with second-order finite-differences.  We test this with both periodic and stationary boundary conditions, initial conditions
\begin{align}
  \label{eq:brusselator2d_ic}
  u(x,y,0) = 22 \left(y-y^2\right)^{1.5}, \quad v(x,y,0) = 27 \left(x-x^2\right)^{1.5},
\end{align}
and parameters $d=\{0.01,0.1\}$, $A=1$, $B=\{3, 30, 300\}$, $\mathbf{c}_u=(-0.5, 1)$, and $\mathbf{c}_v=(0.4, 0.7)$.  We note that the initial condition \eqref{eq:brusselator2d_ic} differs slightly from the one used in \cite{abdullePIROCKSwissknifePartitioned2013a}, since the version used in that paper was continuous but not differentiable under periodic boundary conditions.  Like in \cite{abdullePIROCKSwissknifePartitioned2013a}, we partition this problem so that both the advection and reaction terms are treated explicitly, with the diffusive terms either using an STS method (for PIROCK and ExtSTS) or treated implicitly (for ARK).  We again perform tests using temporal adaptivity with the same tolerances as in the preceding 1D tests.  We compare all results against a reference solution obtained using the fourth order ARK4(3)7L[2]SA$_1$ method from \cite{KenCarp:19}, and with a tighter relative tolerance of $10^{-8}$.

\begin{figure}
  \centering
  \includegraphics[trim={0 20 330 20}, clip, width=0.38\textwidth]{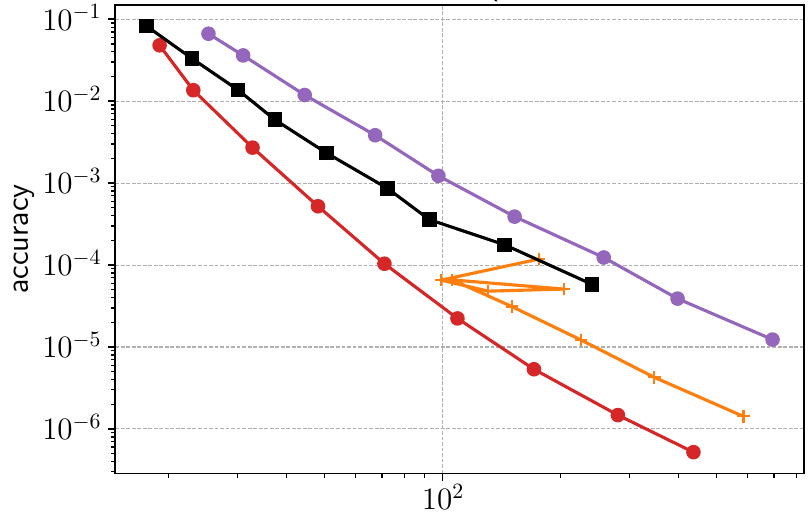}
  \includegraphics[trim={20 20 330 20}, clip, width=0.36\textwidth]{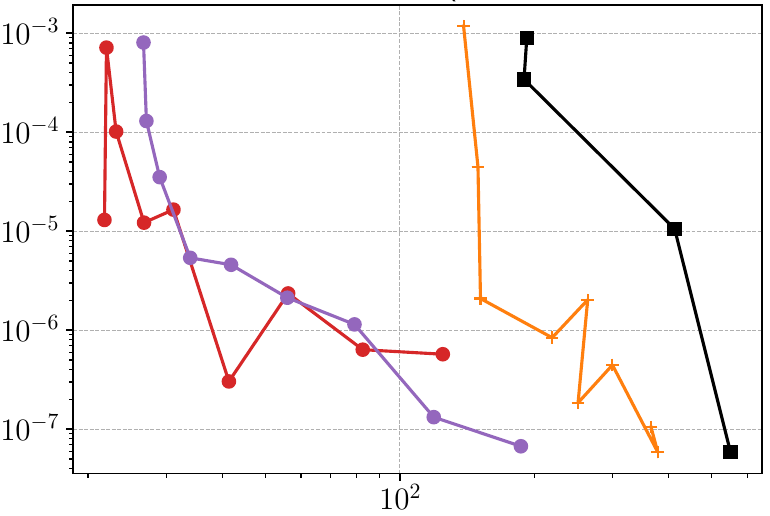}
  \includegraphics[trim={400 30 150 20}, clip, width=0.17\textwidth]{adr2D_stationary_adaptive_runtime_efficiency_d0.01_B3.0e+00.pdf}\\
  \includegraphics[trim={0 20 330 20}, clip, width=0.38\textwidth]{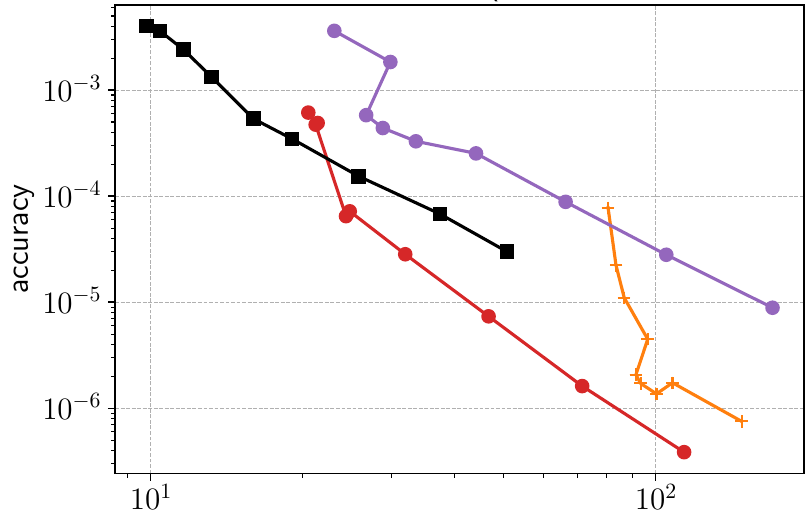}
  \includegraphics[trim={20 20 330 20}, clip, width=0.36\textwidth]{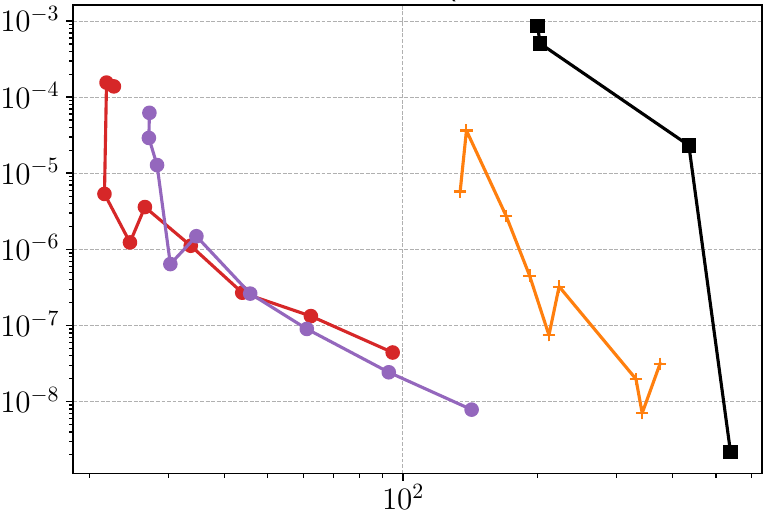}
  \includegraphics[trim={400 0 150 150}, clip, width=0.17\textwidth]{adr2D_stationary_adaptive_runtime_efficiency_d0.01_B3.0e+01.pdf}\\
  \includegraphics[trim={0 0 330 20}, clip, width=0.38\textwidth]{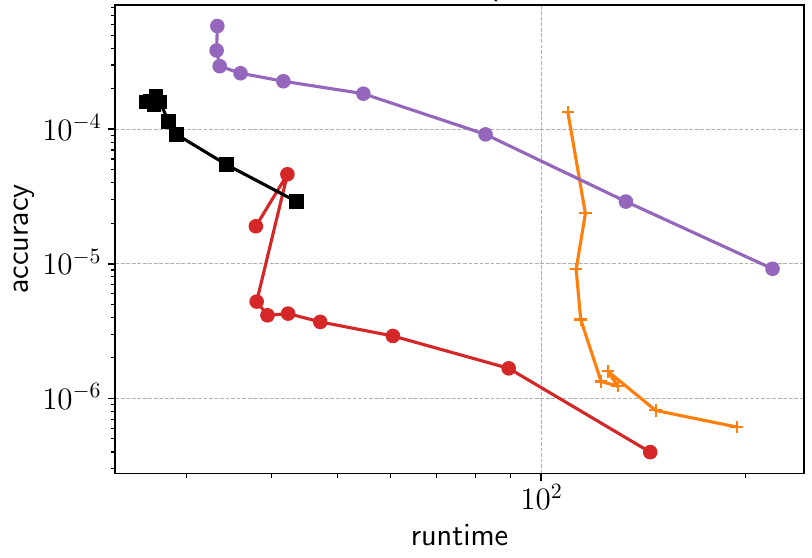}
  \includegraphics[trim={20 0 330 20}, clip, width=0.36\textwidth]{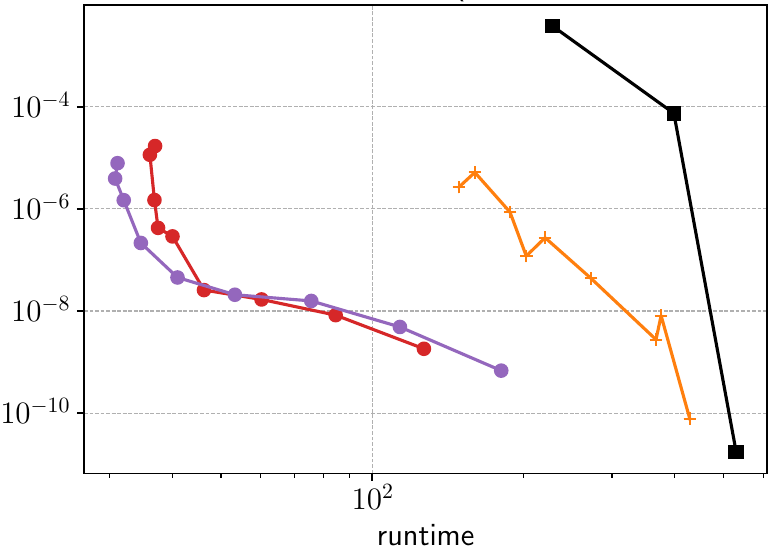}
  \includegraphics[trim={400 0 150 150}, clip, width=0.17\textwidth]{adr2D_stationary_adaptive_runtime_efficiency_d0.01_B3.0e+02.pdf}

  \caption{Runtime efficiency of the adaptive methods on the 2D advection-diffusion-reaction test problem with weaker diffusion parameter $d=0.01$.  The rows show the results for different values of $B$: row one has $B=3$, row two has $B=30$, and row three has $B=300$.  Periodic boundary conditions are shown on the left, and stationary boundary conditions on the right.}
  \label{fig:adr2D_d01}
\end{figure}

In Figure \ref{fig:adr2D_d01} we plot the runtime efficiency of each adaptive method for the problem with weaker diffusion coefficient $d=0.01$, as was used in \cite{abdullePIROCKSwissknifePartitioned2013a}.  For periodic boundary conditions and the least stiff reactions ($B=3$, similar to the case in \cite{abdullePIROCKSwissknifePartitioned2013a}), the Giraldo-ExtSTS method is the most efficient, followed by the PIROCK method at looser tolerances, followed by SSP32-ExtSTS; for this case the ARK method is also quite efficient for tighter tolerances.  As the stiffness increases, however, PIROCK takes over as the most efficient method for loose tolerances, with Giraldo-ExtSTS the best choice at tighter tolerances.  Both are followed by SSP32-ExtSTS, and the ARK method performance degrades by comparison.  For the more challenging case of stationary boundary conditions, we again see poor runtime performance by the PIROCK method, that performs uniformly worse than ARK, which itself is far less efficient than either Giraldo-ExtSTS or SSP32-ExtSTS in performance.  Notably, we see comparable performance in both ExtSTS methods, with SSP32-ExtSTS overtaking Giraldo-ExtSTS for the stiffest case of $B=300$.

\begin{figure}
  \centering
  \includegraphics[trim={0 20 330 20}, clip, width=0.38\textwidth]{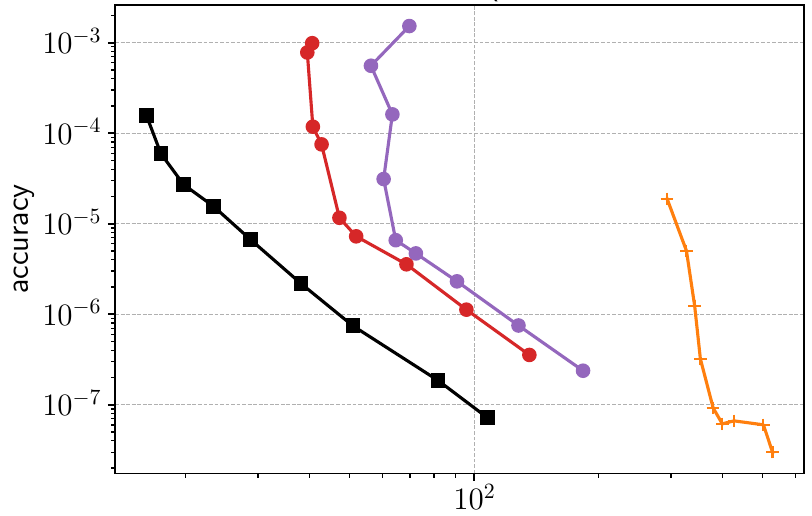}
  \includegraphics[trim={20 20 330 20}, clip, width=0.36\textwidth]{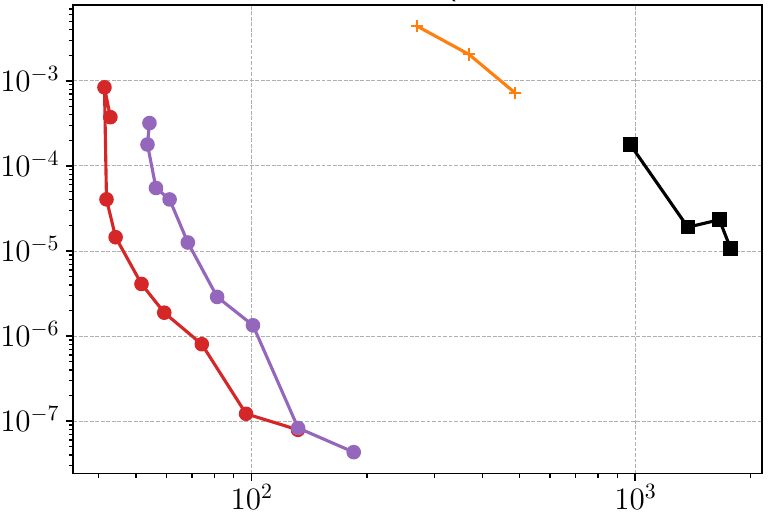}
  \includegraphics[trim={400 30 150 20}, clip, width=0.17\textwidth]{adr2D_stationary_adaptive_runtime_efficiency_d0.10_B3.0e+00.pdf}\\
  \includegraphics[trim={0 20 330 20}, clip, width=0.38\textwidth]{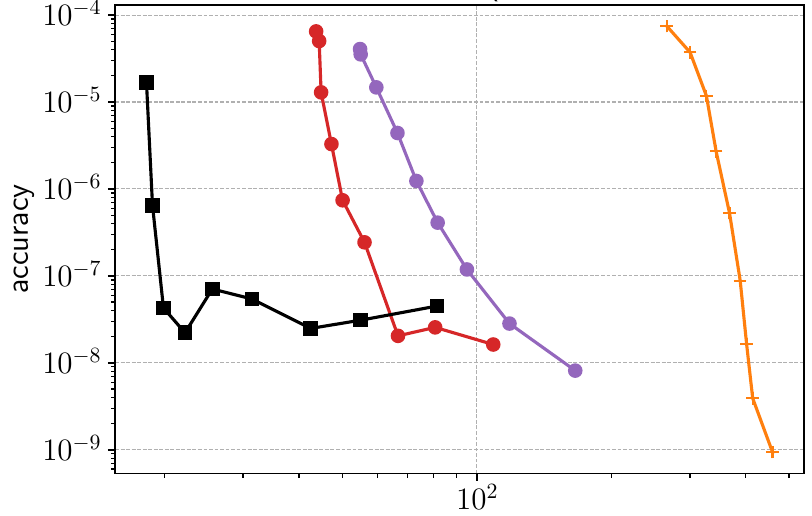}
  \includegraphics[trim={20 20 330 20}, clip, width=0.36\textwidth]{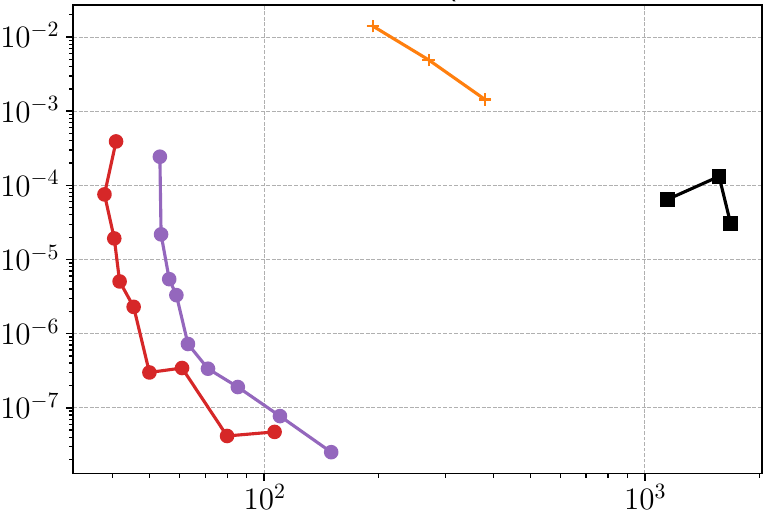}
  \includegraphics[trim={400 0 150 150}, clip, width=0.17\textwidth]{adr2D_stationary_adaptive_runtime_efficiency_d0.10_B3.0e+01.pdf}\\
  \includegraphics[trim={0 0 330 20}, clip, width=0.38\textwidth]{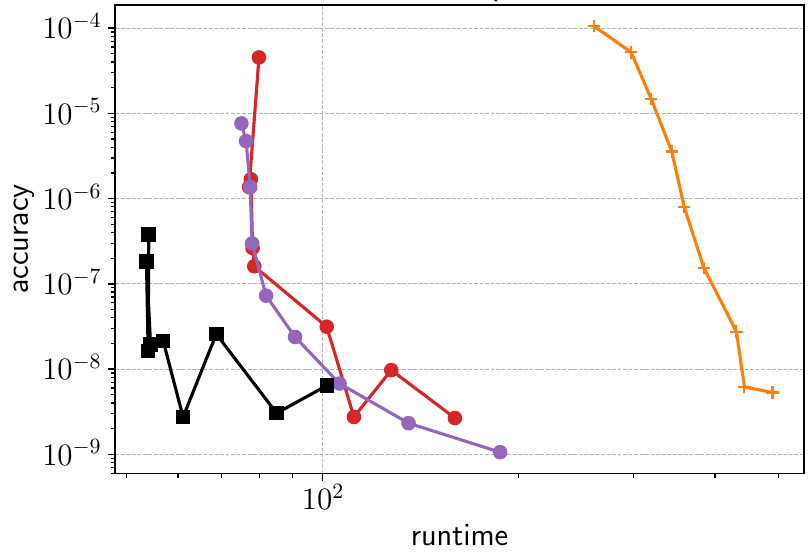}
  \includegraphics[trim={20 0 330 20}, clip, width=0.36\textwidth]{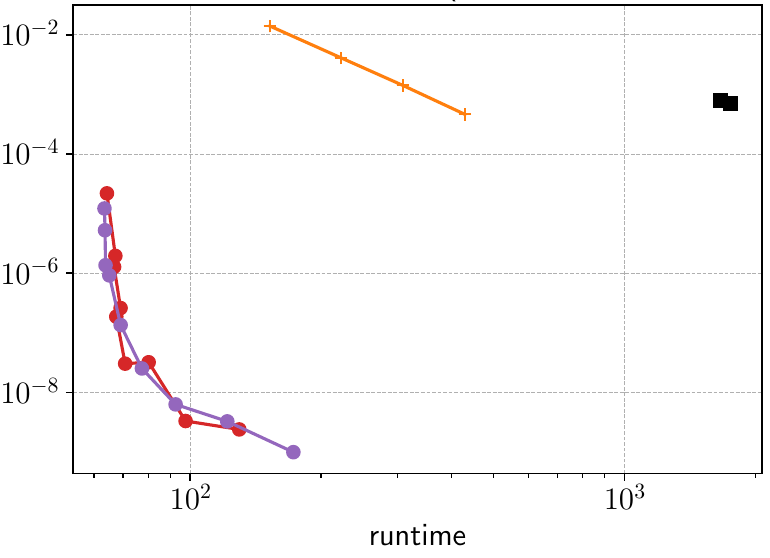}
  \includegraphics[trim={400 0 150 150}, clip, width=0.17\textwidth]{adr2D_stationary_adaptive_runtime_efficiency_d0.10_B3.0e+02.pdf}

  \caption{Runtime efficiency of the adaptive methods on the 2D advection-diffusion-reaction test problem with stronger diffusion parameter $d=0.1$.  The rows show the results for different values of $B$: row one has $B=3$, row two has $B=30$, and row three has $B=300$.  Periodic boundary conditions are shown on the left, and stationary boundary conditions on the right.}
  \label{fig:adr2D_d10}
\end{figure}

In Figure \ref{fig:adr2D_d10} we consider the case with stronger diffusion coefficient $d=0.1$.  Here, for periodic boundary conditions the PIROCK method is clearly the most efficient for all reaction parameters $B$, closely followed by the Giraldo-ExtSTS and then SSP32-ExtSTS methods.  All of the STS-based methods significantly outperform ARK.  For the more challenging stationary boundary conditions, however, PIROCK performs far worse than all other methods, requiring both more runtime and resulting in worse error, than the competing methods.  The ARK methods do not perform much better, struggling to achieve the requested tolerances, and taking more computational effort than the ExtSTS methods, that again perform remarkably well for all values of $B$.  As with the weaker diffusion coefficient, Giraldo-ExtSTS outperforms SSP32-ExtSTS for all but the stiffest case of $B=300$, where they perform similarly.

\section{Conclusions}
\label{sec:conclusions}

We propose a new family of implicit-explicit extensions to super-time-stepping methods, which we call ExtSTS methods, that allow for the efficient solution of advection-diffusion-reaction problems.  These methods apply optimal time integration methods to each component, super time stepping for diffusive components, explicit Runge--Kutta for nonstiff components, and diagonally implicit methods for stiff (but perhaps spatially local) components.  These methods are second-order accurate, include embeddings for robust temporal adaptivity, and achieve significantly smaller errors than existing methods for the same computational effort, due to their tighter coupling between components.

We provide a straightforward algorithm for constructing ExtSTS methods using either existing Runge--Kutta methods or any MRI method, along with any super time stepping sub-method.  Based on this infrastructure, we provide candidate second-order and adaptive implicit-explicit, explicit, and implicit ExtSTS methods.

We demonstrate the performance of these methods on a series of 1D and 2D advection-diffusion-reaction problems, comparing against existing methods including ARK, PIROCK, and Strang splitting.  We find that the proposed ExtSTS methods perform well across a wide range of physical parameters and for both periodic and stationary boundary conditions, in both one and two dimensions.  This is in contrast to the competing methods, which only perform well for some parameter and boundary condition combinations.  Clues to explain these performance differences are visible in the fixed time step tests.  There, for tests with periodic boundary conditions both PIROCK and Strang are competitive and sometimes superior to ExtSTS methods, although the ExtSTS method performance there is still stronger than ARK methods.  However, the tests with stationary boundary conditions require stronger interactions between the advection, diffusion, and reaction processes -- there, only the ExtSTS methods are able to robustly compute solutions at all requested step sizes, while both PIROCK and Strang exhibit instability for some step sizes and for others they resulted in significantly higher solution error.  These larger errors manifest in the adaptive runtime efficiency, where ExtSTS method performance remains largely unchanged across parameters and boundary conditions, while PIROCK performance varies wildly.  While ARK methods are able to achieve accurate solutions for all problem parameters and boundary conditions, their reliance on globally coupled implicit solves renders them less efficient that the ExtSTS methods, an effect that widens as diffusion coefficients increase.

This work may be extended in several directions.  The most obvious extension is to apply the proposed methods to more complex physical problems, such as large-scale problems involving multiple species.  Second, we note that the proposed methods are only second-order accurate, and thus future work may focus on the development of higher-order ExtSTS methods through the construction of higher-order STS methods.  Finally, although the proposed methods may use STS methods with many internal stages, they still evolve all components using a shared time step size.  Thus, future work may focus on the development multi-physics extensions of STS methods that allow multirate subcycling of faster time scales.

\bibliographystyle{siamplain}
\bibliography{references}
\end{document}